\documentclass{article}
\usepackage{amssymb,amsmath,amsthm}
\usepackage{bm}
\usepackage{graphicx}
\usepackage{overpic}

\def\bN{\mathbb N}

\def\bN{{\mathbb N}}

\def\bR{{\mathbb R}}
\def\bS{{\mathbb S}}

\def\cC{{\cal C}}
\def\cD{{\cal D}}
\def\cE{{\cal E}}

\def\cN{{\cal N}}

\def\cS{{\cal S}}

\def\RPt{\bR\hbox{\rm P}^2}
\def\RPT{\bR\hbox{\rm P}^3}

\def\RPn{\bR\hbox{\rm P}^n}

\def\rmU{\uppercase\expandafter{\romannumeral1}}
\def\rmD{\uppercase\expandafter{\romannumeral2}}
\def\rmT{\uppercase\expandafter{\romannumeral3}}

\def\fS{{\mathfrak S}}
\newtheorem{definition}{Definition}
\newtheorem{theorem}{Theorem}
\newtheorem{prop}{Proposition}
\newtheorem{remark}{Remark}

\newtheorem{corollary}{Corollary}
\newtheorem{lemma}{Lemma}
\newtheorem{alg}{Algorithm}
\newtheorem{conj}{Conjecture}

\title{On a generalized Sierpi\'nski fractal in $\RPn$}

\author{Roberto De Leo}

\begin{document}
\maketitle

\begin{abstract}
We associate a fractal in $\RPn$ to each vector basis of $\bR^{n+1}$ and we study 
its measure and asymptotic properties. Then we discuss and study numerically in detail 
the cases $n=1,2,3$, evaluating in particular their Hausdorff dimension.
\end{abstract}

\pagestyle{myheadings}
\thispagestyle{plain}
\markboth{R. De Leo}{A generalized Sierpi\'nski gasket in $\RPn$}
\section{Introduction}
In this paper we study an algorithm that takes a basis of $\bR^{n+1}$ and builds, out of it,
a fractal in $\RPn$. We do this by using the following two basic facts: 
1. In $\RPn$, every $n+1$ points $\{p_i\}$ which do not lie on the same affine hyperplane 
determine a partition of $\RPn$ in the $2^n$ projective $n$-simplices having the points 
$p_i$ as vertices; 2. Given a vector basis $\cE=\{e_i\}$ in $\bR^{n+1}$,
we can build $n+1$ new bases $\cE_i$ by fixing the $i$-th vector $e_i$ and summing it to 
the $n$ remaining ones. Now, consider a vector basis $\cE$ of $\bR^{n+1}$. 
Its vectors $e_i$ projects into $n+1$ points\footnote{Here and throughout 
the paper we denote by $[e]\in\RPn$ the direction of the vector $e\in\bR^{n+1}$.} 
$[e_i]\in\RPn$ and therefore determine a partition of $\RPn$ in projective $n$-simplices
as in point (1). We denote by $S(\cE)$ the one the point $[e_1+e_2+e_3]$. Next, consider 
the $n+1$ projective $n$-simplices $S(\cE_i)$ corresponding to the bases $\cE_i$ defined in (2). 
Their union $\cup_{i=1}^{n+1}S(\cE_i)$ can be thought as the difference between $S(\cE)$ 
and the interior of the projective polytope $Z(\cE)$ (the {\sl body} of $\cE$) having 
the points $[e_i+e_j]$ as vertices. 
By repeating recursively this step on the $S(\cE_i)$, we end up building a $(n+1)$-ary 
tree of bases $T(\cE)=\{\cE_I\}$, where $I=i_1\dots i_k$ is a multiindex, and fractal 
$F(\cE)$ whose points are the ones left inside $S(\cE)$ after removing all the bodies 
$Z(\cE_I)$. 

Topologically, this fractal coincide with the multi-dimensional generalization 
of the Sierpinski triangle~\cite{Sie15}, namely the fractal generated by removing from 
a $n$-simplex $S$ the polytope $Z$ having as vertives the middle points of the edges of $S$. 
Geometrically though they are different because the vertices $[e_i+e_j]$ of the body $Z(\cE)$ 
are closer to the vertices of $S(\cE)$ corresponding to the vectors of higher
Euclidean norm, and even if we start with a basis where all vectors have the
same norm they will not be anymore so after the first step.

We were motivated to study this fractal by the following two reasons. First, this fractal
is invariant with respect to a set of $n+1$ transformations (more specifically, projective 
diffeomorphisms) but they fail to be an Iterated Function System (IFS) because they
are not contractions. In particular the machinery developed for IFSs cannot be applied to 
this case and no analytical bounds are known for the Hausdorff dimension of $F(\cE)$ 
(numerical estimates for the cases $n=2,3$ can be found in Section~\ref{sec:num}). Second, 
the construction leading from $\cE$ to $F(\cE)$ is the $n$-dimensional generalization 
of the structure discovered by the author and I.A. Dynnikov in a fractal, in $\RPt$, 
which describes the asymptotic behaviour of the plane sections of the triply-periodic cubic 
polyhedron
$\cC=\{4,6|4\}$~\cite{DD09}\footnote{We discovered later that this fractal had been already 
considered in the past by G.~Levitt~\cite{Lev93} while studying dynamical systems on the circle.}. 
Indeed one of the results of Novikov's theory of plane sections of triply-periodic 
surfaces~\cite{Nov82,Dyn99} is the following. Consider a connected triply-periodic surface $S$ 
which divides $\bR^3$ in two components which are equal modulo translations. Then there exists
no bundle of parallel planes whose intersections with $S$ are all compact. The asymptotics of
the open (i.e. non-compact) sections, as function of the direction of the bundle of planes, 
are described by a ``labeled cut-out fractal'' in the projective plane (see~\cite{DeL03a,DeL06} 
for a few other concrete examples).
Here by {\sl labeled cut-out fractal} we mean a fractal $F$ which is obtained by removing, 
from an initial region, a sequence of closed sets $\{Z_i\}$ whose interiors is pairwise 
disjoint and such that to each $Z_i$ it is associated an element $b_i$ of some set of ``labels'' 
$B$.
In the particular case of Novikov's theory above, the set of labels is the set of all 
indivisible triples of integers and the geometrical meaning of the fractal 
is the following. Let $d$ be any vector in $\bR^3$ whose direction $[d]$ belongs to some $Z_i$ 
labeled by $b_i$. Then the open sections arising by cutting the surface $S$ with planes 
perpendicular to $d$ are strongly asymptotic to a straight line whose direction is given 
by the ``vector product'' $d\times b$.
Going back to the surface $\cC$, it turns out that its corresponding fractal $F_\cC$ 
essentially coincides with the fractal $F(\cE)\subset\RPt$ associated to the basis 
$\cE=\{(1,0,1),(0,1,1),(1,1,0)\}$ of $\bR^3$.
We conclude this long digression by pointing out that the Novikov's theory of
plane sections of triply-periodic surfaces is the mathematical model for the phenomenon 
of the anisotropic behaviour of magnetoresistance in normal metals at low temperature 
and under a strong magnetic field (see~\cite{LP60} and~\cite{NM03}
for more details about the physics and the dynamics of this phenomenon). In particular
$F_\cC$ encodes the information on the conduction of the electric current in a metal 
having $\cC$ as Fermi Surface.

The paper is organized as follows. In Section~\ref{sec:struct} we define the basic
objects and prove a few elementary facts about them. In Section~\ref{sec:zero} we prove
that all fractals $F(\cE)$ have zero volume with respect to some natural measure on the
projective space and study their asymptotic properties, showing in particular that
they are related with the $n$-bonacci sequences. Finally,
in Section~\ref{sec:num}, we discuss in detail the cases $n=1,2,3$ and present numerical 
results indicating that, unlike the Sierpi\'nski case, the Hausdorff measure of $F(\cE)$ 
may be non-integer even for $n=3$. 
\section{Structure of the fractal}
\label{sec:struct}
Let $\cE=\{e_1,\dots,e_{n+1}\}$ be a vector basis of $\bR^{n+1}$ and let us call
{\sl volume} of $\cE$ the Euclidean volume of the $(n+1)$-simplex of $\bR^{n+1}$ 
naturally associated to it. Every of the $n+1$ sets $\cE_j=\{e^j_1,\dots,e^j_{n+1}\}$, 
$j=1,\cdots,n+1$, defined by 
$$
\begin{cases}
  e^j_i=e_i+e_j,i\neq j\cr
  e^i_i=e_i\cr
\end{cases}
$$
is also a vector basis of $\bR^{n+1}$ and has the same volume as $\cE$. 
Repeating recursively this procedure, we get an inifnite $(n+1)$-ary
ordered rooted tree $T(\cE)=\{\cE_{I_k}\}_{k\in\bN}$ of bases of $\bR^{n+1}$,
all with the same volume, with $\cE$ as root. The multi-index 
$I_k=i_1i_2\cdots i_{k-1}i_k$ describes the steps needed to build the basis 
from the root, namely $\cE_{I_{k-1} i_k}=(\cE_{I_{k-1}})_{i_k}$.

This tree structure corresponds to the limit process for building
a fractal on $\RPn$. Indeed let $\cE=\{e_1,\dots,e_{n+1}\}$ and denote
by $[e_i]\in\RPn$ the direction of the vector $e_i$. To $\cE$ it is naturally 
associated a projective $n$-simplex $S(\cE)$ defined in the following way.
The $n+1$ points $[e_i]\in\RPn$ are the vertices of $2^n$ projective $n$-simplices 
whose interiors are pairwise disjoint and whose union gives the whole $\RPn$;
the point $p=[e_1+\dots+e_{n+1}]$ is not a boundary point for any of them and we
denote by $S(\cE)$ the one which contains $p$. Now consider the bases $\cE_i$
at the first recursion level of $T(\cE)$. The vertices of the projective 
$n$-simplex $S(\cE_i)$ are the points  
$\{[e_1+e_i],\dots,[e_{i-1}+e_i],[e_i],[e_{i+1}+e_i],\dots,[e_{n+1}+e_i]\}$,
i.e. $S(\cE_i)$ is contained inside $S(\cE)$, shares with it the vertex
$[e_i]$ and (part of) all the edges coming out from that point and has in 
common exactly one vertex with each other $S(\cE_j)$, $j\neq i$.

Let now $F_1=\cup_{i=1}^{n+1}S(\cE_i)$. The difference between $S(\cE)$ and
$F_1$ is the interior of the projective polytope $Z(\cE)$ having the 
$n(n+1)/2$ points $[e_i+e_j]$, $i\neq j$, as vertices. We call $Z(\cE)$
the {\sl body} of $\cE$.
More generally, let $F_k=\cup_{|I|=k}S(\cE_I)$ be the $k$-th level of recursion 
of the fractal, with $F_0=S(\cE)$. The set $F_k$ is obtained from $F_{k-1}$ 
by erasing the interiors of the $(n+1)^{k-1}$ bodies $Z_J=Z(\cE_J)$, $|J|=k-1$. 
The fractal $F(\cE)$ is then obtained as the limit $F(\cE)=\cap_{k\in\bN}F_k$.

Note that we can always find an affine $n$-plane (i.e. a canonical chart for 
the projective space) inside $\RPn$ which contains the entire $S(\cE)$ and, 
therefore, the whole $F(\cE)$.
From now on then we will consider often $S(\cE)$ and all the $S(\cE_I)$ as an 
$n$-simplex inside $\bR^n$.
This allows to provide another geometric characterization of the algorithm generating
the fractal. Indeed the $k$-skeleton of $S(\cE)$ is the set of the 
convex hulls associated to the ${n+1\choose k}$ different subsets of $k$ elements of 
$\cE$, namely the convex hulls of the sets $\{[e_{i_1}],\dots,[e_{i_k}]\}$ where no
two indices are equal.
\begin{definition}
  The vector $b(\cE)=\sum_{i=1}^{n+1}e_i\in\bR^{n+1}$ is called the {\sl barycenter} 
  of the $n$-simplex $S(\cE)$, where $\cE=\{e_i\}_{i=1,\cdots,n+1}$.
  Analogously, the barycenter of its $k$-face of vertices $\{[e_{i_1}],\dots,[e_{i_k}]\}$
  is the vector $b_I(\cE)=\sum_{j=1}^{k}e_{i_j}\in\bR^{n+1}$, where $I=i_1\dots i_k$.
  By abuse of notation we sometimes call barycenter its direction
  $[b]\in \RPn$. It will be clear from the context which one we are referring to.
\end{definition}
%
%
\begin{lemma}
  Let $\cE$ be a basis of $\bR^{n+1}$ and $f_{I_k}$ the $k$-subsimplex of $S(\cE)$
  corresponding to $\{e_{i_1},\dots,e_{i_k}\}\subset\cE$.
  Then the projection of the barycenter $[b_{I_k}]$ of $f_{I_k}$ from the vertex 
  $[e_{i_j}]$ on the $(k-1)$-face of $f_{I_k}$ opposite to it, namely the one 
  corresponding to the $k-1$ vectors $\{e_{i_1},\dots,e_{i_k}\}\setminus\{e_{i_j}\}$,
  coincides with the barycenter of that face.
\end{lemma}
\begin{proof}
  This relation is clearly recursive and therefore it is enough to prove the theorem 
  in the case of the $n$-simplex $S(\cE)$ and any of its faces.
  Let us consider what happens for the vertex $[e_1]$: the face $f_1$ of $S(\cE)$ 
  opposite to it corresponds, in $\bR^{n+1}$, to the $n$-plane $\hat f_1$ spanned 
  by the $n$ vectors $\cE^{(1)}=\{e_j\}_{j\neq1}$ and the line $l_1$ joining $[e_1]$ 
  to $[b]$ corresponds to the 2-plane $\hat l_1$ spanned by $e_1$ and $b$. 
  Since $b=\sum_{k=1}^{n+1}e_k$, clearly the only linear combinations belonging to both $\hat f_1$ 
  and $\hat l_1$ are the span of the vector $\sum_{k=2}^{n+1}e_k = b - e_1$. In other
  words, the intersection between $f_1$ and $l_1$ is $[e_2+\dots+e_{n+1}]$, which 
  is indeed the barycenter of $\cE^{(1)}$ and similarly for the other vertices.
\end{proof}
\begin{prop}
  The vertices of the body $Z(\cE)$ corresponding to a basis $\cE$ can be obtained 
  in the following way: project the barycenter of $\cE$ from its vertices to its 
  faces and repeat recursively this procedure until the edges are reached. 
  The $n(n+1)/2$ points obtained are the vertices of $Z(\cE)$.
\end{prop}
\begin{proof}
  The recursive procedure makes sense because, thanks to the previous lemma, we know that
  the projection of the barycenter on a face via the vertex opposite to it coincides 
  with the barycenter of the face. When we reach the edges, therefore, we are left with
  their barycenters, which are clearly the $n(n+1)/2$ points $[e_i+e_j]$, $i\neq j$.
\end{proof}
Finally, we provide a third way to describe this fractal. Recall that the Sierpi\'nski gasket 
and its natural multi-dimensional generalization can be seen the invariant set of a 
Iterated Functions Systems (IFS). Similarly, we prove below that the fractal $F(\cE)$ 
is the invariant set of $n+1$ projective diffeomorphisms $\{\psi_i\}$. They do not form 
however, strictly speaking, a IRS because they are not contractions; in particular the 
Jacobian of each of them is the identity in the omonimous vertex $[e_i]$ of $S(\cE)$.
%
\begin{prop}
  The fractal $F(\cE)$ is invariant with respect to the $(n+1)$ projective automorphisms 
  $\psi_i$ of $\RPn$ induced by the linear transformations $\hat\psi_i$ defined by 
  $\hat\psi_i(e_j)=e_i+e_j$, $j\neq i$, $\hat\psi_i(e_i)=e_i$.
\end{prop}
\begin{proof}
  This is simply a consequence of the fact that the tree itself $T(\cE)$ is clearly invariant 
  under the action of the $\hat\psi_i$, so that the $\psi_i$ map the set of bodies
  $Z(\cE_I)$ into itself and therefore leave the fractal invariant.
\end{proof}
\begin{remark}
  Every body $Z_I(\cE)$, $|I|=k$, is the image of the root body $Z(\cE)$ via the map
  $\psi_I:=\psi_{i_k}\circ\dots\circ\psi_{i_1}$.
\end{remark}
%
%
\section{Measure and Asymptotics of the fractal}
\label{sec:zero}
We start by proving that the volume of $F(\cE)$ is zero with respect to any measure $\mu$ 
induced on $\RPn$ by the Lebesgue measure on any affine $n$-plane or, equivalently, with 
respect to the measure induced by the canonical one on the sphere.

\begin{figure}
  \begin{center}
    \includegraphics[width=11cm]{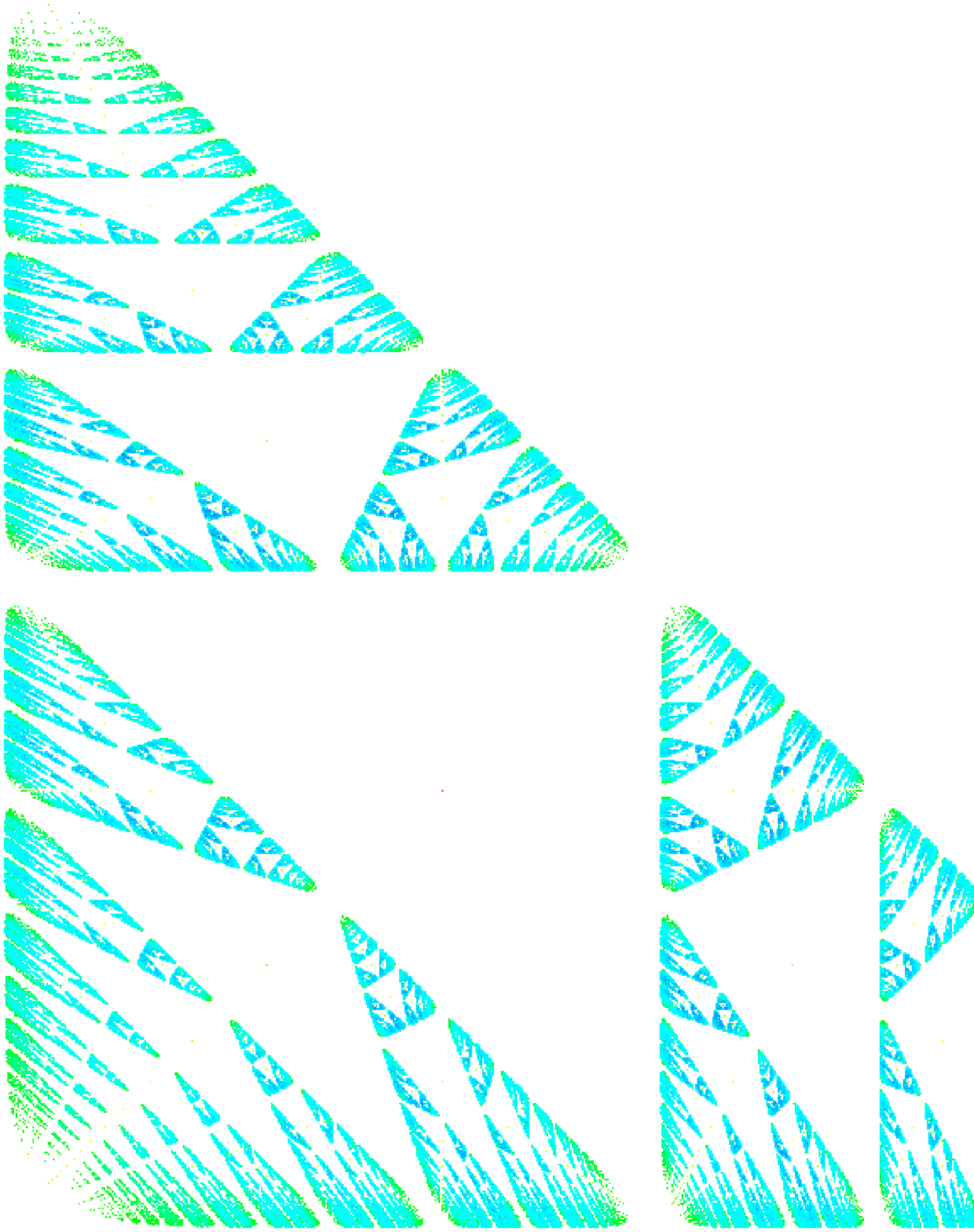}
    \caption{%
      \small
      Plot of the barycenters of all 3487590 bodies up to the thirteenth recursion level
      for $\cE=\{(1,0,1),(0,1,1),(0,0,1)\}$ in the $h^3=1$ projective chart of $\bR \hbox{P}^2$.
      The colors of the points goes from red to blue as the Euclidean norm of the barycenters 
      grows.
      By Corollary~\ref{cor:bardensity} the closure of the set of barycenters is equal to 
      the boundaries of all bodies plus the points of the fractal $F(\cE)$, so this picture
      represents an approximation of the real fractal (and actually no point shown belongs
      to $F(\cE)$ since barycenters are all contained inside the bodies. Nevertheless they
      can approximate as close as wished the set $F(\cE)$ and so they can be used to derive 
      a numerical evaluation of the box-counting dimension of $F(\cE)$.
    }
    \label{fig:labels}
  \end{center}
\end{figure}

We begin with a technical lemma:
\begin{lemma}
\label{tecLemma}
  The maximum $m_{k,n}$ of the functions
  $$
  f_{k,n}(v_1,\cdots,v_n)=\frac{(1+\displaystyle\sum_{i=1}^n v_i)^n}{(1+\displaystyle\sum_{i=1}^n v_i+k(1+\displaystyle\sum_{i=2}^n v_i))(1+\displaystyle\sum_{i=1}^n v_i+(1+k)(1+\displaystyle\sum_{i=2}^n v_i))^n}
  $$
  where $k,n\in\bN$, $n\geq2$, on the $n$-simplex $S$ with vertices in $p_0=(0,\cdots,0)$, $p_1=(1,0,\cdots,0)$,
  $\cdots$, $p_n=(0,\cdots,0,1)$, is given by 
  $$
  m_{k,n}=f_{k,n}(1,0,\cdots,0)=\frac{2^n}{(2+k)(3+k)^n}
  $$ 
  except for the case $n=2$, $k=0$, where
  $
  m_{0,2}=f_{0,2}(0,0,\cdots,0)=1/4
  $ 
\end{lemma}
\begin{proof}
  A direct computation shows that the derivative of $f_{k,n}$ with respect to any $v_i$, $i>1$,
  is negative inside $S$ and therefore the maximum is attained at the smallest values for those
  variables. Then we are left with the function 
  $$h_{k,n}(v_1)=\frac{(1+v_1)^n}{(1+k+v_1)(2+k+v_1)^n}$$
  whose derivative 
  $$h'_{k,n}(v_1)=\frac{(1+v_1)^{n-1}\{k^2 n + (n-2-v_1)(1+v_1) + k[n(2+v_1)-1-v_1]\}}{(1+k+v_1)^2(2+k+v_1)^{n+1} }$$
  in the domain $v_1\in(0,1)$ is always positive for $k>0$ while for $k=0$ is always positive for $n>2$
  and always negative for $n=2$.
  \par
  Hence $f_{k,n}$ will reach its maximum in the origin when $k=0$, $n=2$ and in the point $(1,0,\cdots,0)$ 
  in all other cases.
\end{proof}
%
The following proof is a generalization to any $n$ of the proof provided in~\cite{DD09} for the case $n=2$.
\begin{theorem}
  \label{thm:zeroMeasure}
  The fractal set $F(\cE)$ is a null set for $\mu$. 
\end{theorem}
\begin{proof}
  It is enough to prove the therem for a particular choice of 
  $\cE=\{e_1=(1,0,\cdots,0,1),\cdots,e_n=(0,\cdots,0,1,1),e_{n+1}=(0,\cdots,0,1)\}$. 
  With this $\cE$, the $n$-simplex $S(\cE)$ is contained in the affine plane 
  $\pi=\{h_{n+1}\neq0\}$ with respect to the homogeneous coordinates $[h_1:\dots:h_{n+1}]$. 
  On $\pi$ we use the canonical coordinates $(v_1,\dots,v_n)$ defined by $v_i=h_i/h_{n+1}$
  and the measure 
  $$
  d\mu' = \frac{dv_1\cdots dv_n}{\left(1+\displaystyle\sum_{i=1}^n |v_i|\right)^n}\,.
  $$
  Each of the $\psi_i$ maps the whole fractal $F=F(\cE)\subset S(\cE)$ into disjoint
  sets $\psi_i(F)=F(\cE_i)\subset S_i$ so that
  $\hat\mu(F)=\sum_{i=0}^{n+1}\hat\mu(\psi_i(F))$.
  With our choice of $\cE$ and $\hat\mu$ both $S(\cE)$ and the measure (at least close
  to $S(\cE)$) are invariant with respect to every permutation of the first $n$ basis 
  vectors and therefore 
  $\hat\mu(\psi_i(F))=\hat\mu(\psi_1(F))$, $\forall i=1,\cdots,n$, so that
  $\hat\mu(F)=n\hat\mu(\psi_1(F))+\hat\mu(\psi_{n+1}(F))$.
  By repeating this procedure on $\psi_{n+1}(F)$ we find that 
  $$\hat\mu(\psi_{n+1}(F))=n\hat\mu(\psi_{n+1}\circ\psi_1(F))+\hat\mu(\psi_{n+1}\circ\psi_{n+1}(F))$$
  so that $\hat\mu(F)=n\hat\mu(\psi_1(F))+n\hat\mu(\psi_{n+1}\circ\psi_1(F))+\hat\mu(\psi_{n+1}\circ\psi_{n+1}(F))$ 
  and finally, by recursion,
  $$
  \hat\mu(F)=n\sum_{k=0}^\infty\hat\mu(\psi_{n+1}^k\circ\psi_1(F))
  $$ 
  since $\lim_{k\to\infty}\hat\mu(\psi_{n+1}^k(F))=0$.
  \par
  We will now show that $\hat\mu(\psi_{n+1}^k\circ\psi_1(F)) \leq c^{(n)}_k\hat\mu(F)$ with 
  $\sum_{k=0}^\infty c^{(n)}_k<1/n$, which leads immediately $\hat\mu(F)=0$. 
  Note indeed that, with this particular choice of the basis, the action of the $\psi_i$ 
  on the corresponding homogeneous coordinates is given by
  $$
  \begin{cases}
    \psi_1([h_1:h_2:\cdots:h_n:h_{n+1}])=[h_{n+1}:h_2:\cdots:h_n:2h_{n+1}-h_1]\cr
    \vdots\cr
    \psi_n([h_1:h_2:\cdots:h_n:h_{n+1}])=[h_1:h_2:\cdots:h_{n+1}:2h_{n+1}-h_n]\cr
    \psi_{n+1}([h_1:h_2:\cdots:h_n:h_{n+1}])=[h_1:h_2:\cdots:h_n:\sum_{i=1}^{n+1}h_i]\cr
  \end{cases}
  $$
  Since the fractal is invariant with respect to the projective transformation
  $$
  R([h_1:h_2:\dots:h_n:h_{n+1}]) = [2h_{n+1}-\sum_{i=1}^{n+1}h_i:h_2:\cdots:h_n:h_{n+1}]
  $$
  corresponding to the exchange of the vectors $e_1$ and $e_{n+1}$, 
  we can replace $\psi_{n+1}^k\circ\psi_1(F)$ with $\psi_{n+1}^k\circ\psi_1\circ R(F)$.
  \par
  Then
  $$
  \psi_1\big( R([h_1:h_2:\cdots:h_n:h_{n+1}])\big) = [h_{n+1}:h_2:\cdots:h_n:\sum_{i=1}^{n+1}h_i]
  $$
  and finally
  $$
  \psi_{n+1}^k\Big(\psi_1\big( R([h_1:h_2:\cdots:h_n:h_{n+1}])\big)\Big) = 
  [h_{n+1}:h_2:\cdots:h_n:\sum_{i=1}^{n+1}h_i+k\sum_{i=2}^{n+1}h_i]\,.
  $$
  In the chart $v_i=h_i/h_{n+1}$ the map $f_k=\psi_{n+1}^k\circ\psi_1\circ R$ is represented by 
  $$
  f_k(v_1,\cdots,v_n) = (1/D,v_2/D, \cdots, v_n/D), D=1+\sum_{i=1}^n v_i+k(1+\sum_{i=2}^n v_i)\,.
  $$
  A direct computation shows that the Jacobian of $f_k$ is given by 
  $$\bigg|\det(\frac{\partial f_k^i}{\partial v_j})\bigg| = \frac{1}{D^{n+1}}$$
  so that $\hat\mu(F^{(k)}_1)\leq c^{(n)}_{k-1}\hat\mu(F)$ for
  $$
  c^{(n)}_k = \max_{(v_i)\in S(\cE)}\bigg|\det(\frac{\partial f_k^i}{\partial v_j})\bigg| \frac{(1+\sum_{i=1}^n v_i)}{(1+\sum_{i=1}^n f_k^i(v_i))}
  = \max_{(v_i)\in S(\cE)}\frac{(1+\sum_{i=1}^n v_i)^{n+1}}{D\cdot(v_1+(k+2)\sum_{i=2}^n v_i)^{n+1}}\,.
  $$
  As shown in Lemma~\ref{tecLemma} $$c^{(n)}_k=\frac{2^n}{(2+k)(3+k)^n}$$ with the sole exception of the case 
  $n=2$, $k=0$, in which case $c^{(2)}_0=1/4$. If $n=2$ then, as already shown in~\cite{DD09}, 
  $$
  \sum_{k=0}^\infty c^{(2)}_k = 1/4 + \sum_{k=1}^\infty \frac{2^{2}}{(2+k)(3+k)^{2}} = \frac{253}{36}-\frac{2}{3}\pi^2\simeq0.45<\frac{1}{2}\,.
  $$  
  In the $n>2$ case instead we use the fact that
  $$
  \vbox{
    \halign{\hfill$\displaystyle #$&$\displaystyle #$\hfill\cr
      \sum_{k=0}^\infty c^{(n)}_k =& \sum_{k=0}^\infty \frac{2^n}{(2+k)(3+k)^n} < 
      \frac{2^{n-1}}{3^n} + 2^n\int_0^\infty\frac{dx}{(2+x)(3+x)^n} =\cr
      =&\frac{2^{n-1}}{3^n} + 2^n\left[\int_0^\infty(\frac{1}{2+x}-\frac{1}{3+x})dx
        -\sum_{k=2}^n\int_0^\infty\frac{dx}{(3+x)^k}\right]=\cr
      =&\frac{2^{n-1}}{3^n} + 2^n\left[\ln\frac{3}{2}-\sum_{k=1}^{n-1}\frac{1}{k 3^k}\right]\,.\cr
    }
  }
  $$
  By Taylor's expansion theorem applied to $\log(1-x)$ we know that there exist a $\xi\in(0,1/3)$ such that
  $$
  \ln\frac{3}{2}=\sum_{k=1}^n\frac{1}{k 3^k}+\frac{1}{(n+1)3^{n+1}(1-\xi)^{n+1}}<\sum_{k=1}^n\frac{1}{k 3^k}+\frac{1}{(n+1)2^{n+1}}
  $$
  so that finally
  $$
  \sum_{k=0}^\infty c^{(n)}_k < \frac{2^n}{3^n}(\frac{1}{2}+\frac{1}{n})+\frac{1}{2(n+1)}\,.
  $$
  It is easy to verify that the analytical function $g(x)=(\frac{2}{3})^x(\frac{1}{2}+\frac{1}{x})+\frac{1}{2(x+1)}$
  is bigger than $h(x)=1/x$ for $x\geq4$, which proves that $\sum_{k=0}^\infty c^{(n)}_k<1/n$ for all $n\geq4$.
  \par
  We complete the proof by verifying the case $n=3$ by a direct computation:
  $$
  \sum_{k=0}^\infty c^{(3)}_k = \sum_{k=0}^\infty \frac{2^3}{(2+k)(3+k)^3} = 13 - \frac{4}{3} \pi^2 - 8 \zeta(3) \simeq 0.22 < \frac{1}{3}\,.
  $$
\end{proof}
%
Next Corollary will be used later to justify one of the numerical methods we used to evaluate the
box-counting dimension of the fractal. An illustration of it can be found in Fig.~\ref{fig:labels}.
\begin{corollary}
  \label{cor:bardensity}
  The fractal set $F(\cE)$ is contained in the set of accumulation points of the set of barycenters.
  In particular, the closure of the set of the barycenters is equal to the union of $F(\cE)$ with
  the boundaries of the bodies $Z(\cE_I)$, $\cE_I\in T(\cE)$.
\end{corollary}
\begin{proof}
  Since $F(\cE)$ has zero measure it cannot contain any open set. In other words, every open set 
  inside $S(\cE)$ either is contained inside a body $Z(\cE_I)$ for some multi-index $I$ or contains 
  one of them. Let $p\in F(\cE)$. Then any open neighborhood of $p$ is not contained inside a body
  and therefore contains one. Inside every body lies a barycenter and so $F(\cE)$ is contained 
  in the closure of the (countable) set of barycenters.
\end{proof}
%
%
%
In order to study the asymptotics of the fractal it is convenient to pose the following 
definition:
\begin{definition}
  We call ``section'' of an infinite tree $T$ a sequence $\{t_i\}_{i\in\bN}\subset T$ such that 
  each element $t_n$ (except for the first) is child of its antecedent $t_{n-1}$.
\end{definition}
In the Sierpinski case the asymptotics properties do not depend on the particular section
but in case of $F(\cE)$ they do. E.g. consider an edge $t$ of $T(\cE)$, i.e. a section 
$t=\{t_k\}$ defined by $t_k=\cE_{\underbrace{ii\cdots i}_k}$ for some index $i$. Then the volume 
of the simplices $S(t_k)$ decreases polynomially with $k$, while in the Sierpinski case
they always decrease exponentially. Below we study the sections
where the volumes grow faster. They are related to $n$-bonacci sequence, namely sequences
whose $k$-th element is equal to the sum of the previous $n$ ones, thanks to the next 
proposition:
\begin{prop}
  \label{thm:bar}
  Let $A$ be the set of bodies having non-empty intersection with $Z_I(\cE)$ and $B$
  the set of indices of the hyperfaces of $S(\cE)$ (labeled after the index of the vertex
  opposite to it) having points in common with $Z_I(\cE)$.
  Then $b(Z_I(\cE))=\sum_{z\in A}b(z)+\sum_{k\in B}b_k$, where $b_k=\sum_{i=1}^{n+1}e_i - n e_k$. 
\end{prop}
\begin{proof}
  Since the fractal is invariant under the $\psi_i$ and they are induced by linear transformations, 
  it is enough to prove this property for the barycenter $b$ of the root cut-out polytope $Z=Z(\cE)$,
  which cuts all faces of $S=S(\cE)$.
  Using the $\psi_i^{-1}$ it is easy to determine that the body corresponding to the $k$-th face of 
  $S$ is the body of the simplex of vertices 
  $$
  \{[e_1-e_k],\dots,[e_{k-1}-e_k],[e_k],[e_{k+1}-e_k],\dots,[e_{n+1}-e_k]\}
  $$
  and therefore its barycenter is the vector $b_k=\sum_{i=1}^{n+1}e_i - n e_k$.
  Now it is easy to verify that
  $$
  \sum_{k=1}^{n+1}b_k=\sum_{k=1}^{n+1}\left(\sum_{i=1}^{n+1}e_i - n e_k\right) = (n+1)\sum_{i=1}^{n+1}e_i - n\sum_{k=1}^{n+1}e_k 
  = \sum_{i=1}^{n+1}e_i = b\,.
  $$
\end{proof}
\begin{figure}
  \label{fig:spiral}
  \begin{center}
    \includegraphics[width=7cm]{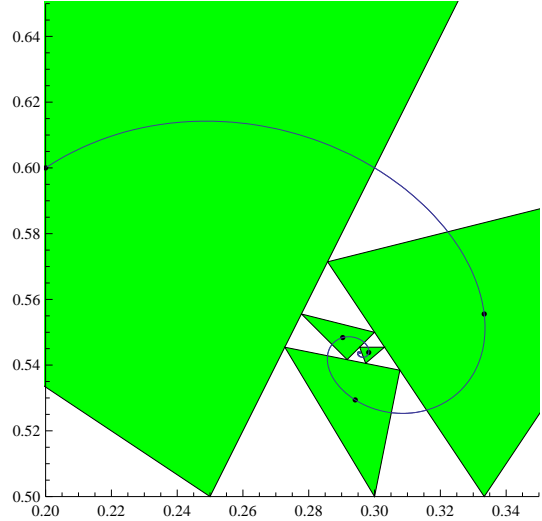}
    \caption{%
      \small
      Detail, in the $h_3=1$ chart, of the first few bodies corresponding to a Tribonacci 
      section starting from the root of the tree $T(\cE)$, where $\cE=\{(1,0,1),(0,1,1),(0,0,1)\}$. 
      The first barycenter $b_1=(1,1,3)$ of the section is not shown. 
      The next five ones, whose projection on $\bR P^2$ is shown above, are
      $b_{12}=(1,3,5)$, $b_{123}=(3,5,9)$, $b_{1231}=(5,9,17)$, $b_{12312}=(9,17,31)$ and
      $b_{123123}=(17,31,57)$.
      The centers of the bodies of the section lie on a smooth ``Tribonacci projective spiral'' 
      drawn above which is winding about $(1/\alpha_3,1/\alpha_3^2)\simeq(.296,.544)$.
    }
  \end{center}
\end{figure}
This result suggests the following interesting way of building sections of a 
tree $T(\cE)$ for a basis $\cE$ of $\bR^{n+1}$. Pick any element $t_1=\cE_{i_1}$ and continue 
the section recursively by taking $t^\pm_{j}=\cE_{i1,i1\pm1,\cdots,i1\pm j}$, where all indices 
are meant modulo $(n+1)$. By construction, the body of $t^\pm_{n+2}$ touches the bodies of all 
of the previouselements of $t^\pm_i$ and therefore its barycenter is given exactly by the sum 
of the barycenters of their bodies, and the same happens for all remaining terms $t^\pm_i$, $i>n+2$.
We call {\sl Fibonacci sections} this particular kind of sections because the sequence of
the corresponding barycenters is a $n$-bonacci sequence. Fibonacci sections of $T(\cE)$ are 
relevant for two reasons: 1. they represent the sections with faster growth hich barycenters 
norms grow faster; 2. they provide a way to get explicit expressions for points in $F(\cE)$.
\begin{theorem}
  \label{thm:fibsec}
  Barycenters of $t_k$'s bodies in a Fibonacci sequence grow in norm as $\alpha^k$, where $\alpha$ is 
  the $(n+1)$-bonacci number (i.e. the highest module root of the equation $x^{n+1}=x^n+\cdots+x+1$). 
  This is the highest growth rate for barycenters' norms on a section of $T(\cE)$.
\end{theorem}
\begin{proof}
  We can assume without loss of generality that $e_i$ is the canonical basis for $\bR^{n+1}$, since asymptotics will 
  not change under the action of a single invertible linear transformation, and we can prove
  the result using the norm $\|v\|_1=\sum_{i=1}^{n+1}|v_i|$ because in finite dimension all
  norms are equivalent. 

  It is well known that the $k$-th term, $k>n+1$, of a $(n+1)$-bonacci sequence can be expressed
  as a linear combination with constant coefficients of the $k$-th powers of the $n+1$ complex 
  roots of the $(n+1)$-bonacci equation $x^n=x^{n-1}+\cdots+x+1$. The highest module root is
  known to be real and it is called $(n+1)$-bonacci constant. Asymptotically only the highest
  module root is relevant and this proves the first part of the theorem.
  
  Now, assume that up to the $n$-th recursive step it happens that at
  each step $k$ the bodies with higher baricentric norm are the ones built starting from $Z(\cE)$
  and belonging to a Fibonacci section: then at the following recursive step the bodies with higher
  barycentric norm are exactly the ones which continue those Fibonacci sections.
  Indeed no body can touch more than one body from each tree level since bodies corresponding to the
  same level belong to distinct simplices; hence at the $(k+1)$-th level the bodies' barycenters 
  of the members of those Fibonacci sections are obtained by summing of the highest norm barycenters 
  and the components are all positive, so their norm is the biggest achievable.
\end{proof}
\begin{remark}
  Proposition~\ref{cor:bardensity}, applied to Fibonacci sections, grants that the limit point 
  of a Fibonacci section must belong to $F(\cE)$. Consider for example the Fibonacci sequence 
  generated by 
  $$b_{-n}=(1,\cdots,1,1-n),\cdots, b_{0}=(1-n,1,\cdots,1), b_{1}=(1,\cdots,1)$$ 
  In this case all components follow the very 
  same sequence but the component $j$ is shifted by one with respect to the component $j+1$
  for $j=1,\cdots,n$, namely $b_k^i=b_{k-1}^{j+1}$. The last component $b_k^{n+1}$ has ``initial 
  conditions'' $b_{-n}^{n+1}=-n$, $b_{-n+1}^{n+1}=1$, $\cdots$, $b_{0}^{n+1}=1$, so that
  the first terms of the sequence are $b_1^{n+1}=1$, $b_2^{n+1}=n+1$, $b_3^{n+1}=2n+1$
  and so on. 
  Since the $k$-th term of a $n$-bonacci sequence behaves asymptotically like $\alpha^k$, 
  in $\RPn$ the sequence of the corresponding points converges to $(1:\alpha:\cdots:\alpha^n)$.
\end{remark}
The following theorems shows that edges and Fibonacci sections are respectively the slower
and faster sections with respect to volumes' growth.
\begin{theorem}
  Let $\cE$ be a basis of $\bR^{n+1}$, $T(\cE)=\{t_I\}$ its tree of bases and 
  $B(\cE)=\{b_I\}$ the corresponding tree of barycenters. Then
  there exist real constants $A$, $B$ such that
  $$
  A|I| \leq \|b_I\| \leq B\alpha^{|I|}
  $$
  for all multiindices $I$.
\end{theorem}
\begin{proof}
  We can prove without loss of generality the theorem by fixing the basis as the canonical basis of 
  $\bR^{n+1}$ and the norm as the maximum norm $\|v\|_\infty=\max|v_i|$.

  As shown in Theorem~\ref{thm:fibsec}, the biggest barycenters at every level $k$ are those belonging
  to a Fibonacci section starting by the root element of the tree; the explicit expression for those 
  sections, modulo permutations, is $b_k=(a_{k-n-1},\cdots,a_k)$, $k>n+1$, where 
  $a_k=\sum_{i=1}^{n+1}\lambda_i \alpha_i^k$, the $\alpha_i$ are the root of the $(n+1)$-bonacci equation 
  and $\lambda_i=1/\Pi_{j\neq i}(\alpha_i-\alpha_j)$. 
  We order the roots so that $\alpha_1=\alpha$ is the $(n+1)$-bonacci constant. Hence, for $k$ big enough,
  $$
  \|b_k\|_\infty=|a_k|\leq 2\lambda_1 \alpha^k
  $$
  
  The slowest growth, again modulo permutations, is obtained by those $n$-simplices corresponding 
  to the bases $\{e_1,e_2+ke_1,\cdots,e_{n+1}+ke_1\}$, whose barycenter
  $b_k=nke_1+\sum_{i=1}^{n+1}e_i$ has norm $\|b_k\|_\infty=nk+1$.
\end{proof}
%
Now we provide bounds for the volumes of the bodies $Z(\cE_I)$ in terms of the norms of the
barycenters.
\begin{lemma}
  \label{tecLemma2}
  Let $W=(w_1,\cdots,w_{n+1})\in\bR^{n+1}$, $n>1$, be a vector with non-negative components and 
  let us build out of it a tree $T(W)$ using the same algorithm used to build $T(\cE)$, so that e.g. 
  at the first tree level we find $W_1=(w_1,w_1+w_2,\cdots,w_1+w_{n+1})$ and the other $n$ 
  vectors obtained similarly. Then if the components of $W$ satisfy the inequalities
  \begin{equation}
    \label{eq:ineq}
    \sum_{j\neq j_1,j_2}w_{j}\leq (n-1)(w_{j_1}+w_{j_2})
  \end{equation}
  the same inequalities hold for all other vectors of the tree. 
\end{lemma}
\begin{proof}
  We prove the lemma by induction. 
  Let us assume that the inequality is valid for all vectors
  up to the $k$-th tree level and be $W'=(w'_1,\cdots,w'_{n+1})$ one of the vectors at the level $k$.
  For the symmetry of the problem it is enough to verify that the inequality remains true
  for its first child $W''=W'_1$ and it is enough to check it in any two cases when its first component
  $w''_1$ appears on the right side of the inequality and when it does not.\par

  In the first case let us assume $j_1=1$ and $j_2=2$. Then the inequality reads
  $$\sum_{j=3}^{n+1}w''_j\leq(n-1)(w''_1+w''_2)$$
  that is equivalent to
  $(n-1)w'_1+\sum_{j=3}^{n+1}w'_j\leq(n-1)(2w'_1+w'_2)$
  and therefore to
  $\sum_{j=3}^{n+1}w'_j\leq(n-1)(w'_1+w'_2)$
  which holds by the inductive hypothesis.
  
  In the second case let us assume $j_1=2$ and $j_2=3$. Then the inequality reads
  $$w''_1+\sum_{j=4}^{n+1}w''_j\leq(n-1)(w''_2+w''_3)$$
  that is equivalent to
  $(n-1)w'_1+\sum_{j=4}^{n+1}w'_j\leq(n-1)(2w'_1+w'_2+w'_3)$
  and therefore to
  $w'_1+\sum_{j=4}^{n+1}w'_j\leq(n-1)(w'_1+w'_2+w'_3)+w'_1$
  which is true because, by the induction hypothesis, 
  $$w'_1+\sum_{j=4}^{n+1}w'_j\leq(n-1)(w'_2+w'_3)\leq(n-1)(w'_1+w'_2+w'_3)+w'_1\,.$$
  All remaining inequalities are obtained by permuting the indices.
\end{proof}
\begin{theorem}
  \label{thm:bodyVol}
  For every basis $\cE$ of $\bR^{n+1}$ 
  there exist real constants $A$, $B$
  such that
  $$\frac{A}{\|b_I\|^{n+1}} \leq \mu(Z_I) \leq \frac{B}{\|b_I\|^{n+1}}$$
  for almost all multi-indices $I$.
\end{theorem}
\begin{proof}
  For this proof's sake it is convenient to use the same base $\cE$ and measure $\mu$ 
  of Theorem~\ref{thm:zeroMeasure} and  the maximum norm for the barycenters.

  These choices have some important advantages: 
  1. if we call $e_i^j$, $j=1,\dots,n+1$, the components of the vectors $\{e_i\}=\cE$ with 
  respect to the canonical basis of $\bR^{n+1}$, then the vector $W=(e_1^{n+1},\dots,e_{n+1}^{n+1})$ 
  built with the $(n+1)$-th coordinates of the basis vectors changes, when passing from 
  the basis $\cE_I$ to $\cE_{I,i_{k+1}}$, 
  with the same rule illustrated in the Lemma above and satisfies the set of 
  inequalities~(\ref{eq:ineq}); 2. if $[h_1:\cdots:h_{n+1}]$ are the canonical homogeneous 
  coordinates for $\RPn$, $S(\cE)$ is entirely contained in the open set $h_{n+1}\neq0$;
  3. the component $e_i^{n+1}$ is not smaller than any other component for every $i=1,\dots,n+1$;
  4. the expressions for the volume of $S(\cE)$ and $Z(\cE)$ are particularly simple.
  
  
  Now let $(x^1_i,\cdots,x^{n+1}_i)$ be the components
  of the vectors of the basis $\cE_I$, so that the homogeneous coordinates of the $(n+1)$ vertices 
  of the $n$-simplex $S_I$ will be $A_i=[x^1_i:\cdots:x^{n+1}_i]$ and those of its 
  body $Z_I$ will be $B_{ij}=[x^1_i+x^1_j:\cdots:x^{n+1}_i+x^{n+1}_j]$.
  A direct computation shows that
  $$
  \hat\mu(S_I)=\frac{1}{n!\displaystyle\prod_{i=1}^{n+1}x_i^{n+1}}
  $$
  and 
  $$
  \hat\mu(Z_I)=\frac{1}{n!}\sum_{S\in\fS_n}\frac{1}{\displaystyle\prod_{x_i\neq x_j\in S}(x_i^{n+1}+x_j^{n+1})}
  $$
  where $\fS_n$ is the subdivision of $Z_I$ in $n+1$ simplices $S_I^{(k)}$, where each $S_I^{(k)}$ has the 
  same vertices of $S_i$ except for the $k$-th vertex, which is replaced by the barycenter of $S_I$.
  
  Let us consider now one of the simplices $s\in\fS_n$ and let $[e_{i_{j,1}}+e_{i_{j,2}}]$, 
  $j=1,\dots,n+1$, be its vertices. Note that,since all components of the basis vectors are 
  positive and no component is bigger than the last one,  the barycenter's norm is
  $$
  \|b_I\|_\infty=\sum_{l=1}^{n+1}x_l^{n+1}\,.
  $$
  Hence 
  $$
  1\leq \frac{\|b_I\|^{n+1}_\infty}{\Pi_{j=1}^{n+1}(x_{i_{j,1}}^{n+1}+x_{i_{j,2}}^{n+1})}\leq
  \Pi_{j=1}^{n+1}\bigg(1+\frac{\sum_{j\neq j_1,j_2}x_j^{n+1}}{x_{i_{j,1}}^{n+1}+x_{i_{j,2}}^{n+1}}\bigg)
  \leq n^{n+1}\,.
  $$
  Since we never used in our calculation the particular choice of the indices for the simplex $s$,
  these bounds are valid for all of them and therefore
  $$
  \frac{n+1}{n!\|b_I\|^{n+1}_\infty}\leq\hat\mu(Z_I)\leq \frac{(n+1)n^{n+1}}{n!\|b_I\|^{n+1}_\infty}\,.
  $$
\end{proof}
Note that the inequality above does not hold for the $n$-simplices $S_I$: for example, 
in the basis $\cE$ used above the simplices corresponding to the bases 
$$
\cE_k=\{e_1+ke_{n+1},\cdots,e_n+ke_{n+1},e_{n+1}\}
$$
have barycenter $b_k=(1,\cdots,1,nk+1)$ and volume 
$$\mu_k=\frac{1}{n! x_1^{n+1}\cdots x_{n+1}^{n+1}}=\frac{1}{n!(k+1)^n}$$
which therefore is asymptotic to $1/\|b_k\|^n$ rather than to $1/\|b_k\|^{n+1}$.

Numerical and analytical facts suggest that bodies' diameters are bound by the inverse of their
barycenters' norm; in particular it is known to be true for $n=2$ thanks to an indirect proof 
(see Section~\ref{sec:n2}) and it is confirmed by numerical exploration of the $n=3$ case 
(see Section~\ref{sec:n3}).
We are led therefore to the following conjecture:
\begin{conj}
  \label{conj:Diam}
  For every basis $\cE$ of $\bR^{n+1}$ there exist constants $A$ and $B$ such that 
  $$\frac{A}{\|b(Z)\|^\frac{n+1}{n}} \leq |Z| \leq \frac{B}{\|b(Z)\|}$$
  where $Z$ is any body associated to the tree $T(\cE)$ and $|Z|$ its diameter with
  respect to the canonical distance $d([x],[y])=\arccos\frac{<x,y>}{\|x\|\|y\|}$.
\end{conj}
As for the fractal dimension of $F(\cE)$, we could not find any way to evaluate exact non-trivial 
bounds for it; in next section we present the numerical evaluation of it for the cases $n=2,3$.
%
%
\section{Analysis of the cases $n=1,2,3$}
\label{sec:num}
\subsection{The case n=1}
The construction we discussed above does not strictly speaking apply to the $n=1$ case.
E.g. bodies here are simply single points, Theorem ~\ref{thm:bar} does not apply and 
all asymptotics about the measures of bodies have no meaning here.
Nevertheless a few things survive: the tree $T(\cE)$ and its Fibonacci sections can still be built 
and we can study the asymptotics of the lengths of the 1-simplices.

To begin, let us choose 
$$\cE=\{e_1=(1,0),e_2=(1,1)\}$$
The set $F(\cE)$ is invariant with respect to the projective transformations
$$\psi_1([h_1:h_2])=[h_1+h_2:h_1],\psi_2([h_1:h_2])=[2h_1-h_2:h_1]$$
and it is obtained from the segment $[0,1]$ (in the projective chart $x=1$) by
removing a countable set of infinite (rational) points, so that 
it has full measure and therefore $\dim_H F(\cE)=1$.

%
\begin{figure}
  \begin{center}
    \begin{tabular}{cc}
      \begin{overpic}[tics=5,width=5.5cm]{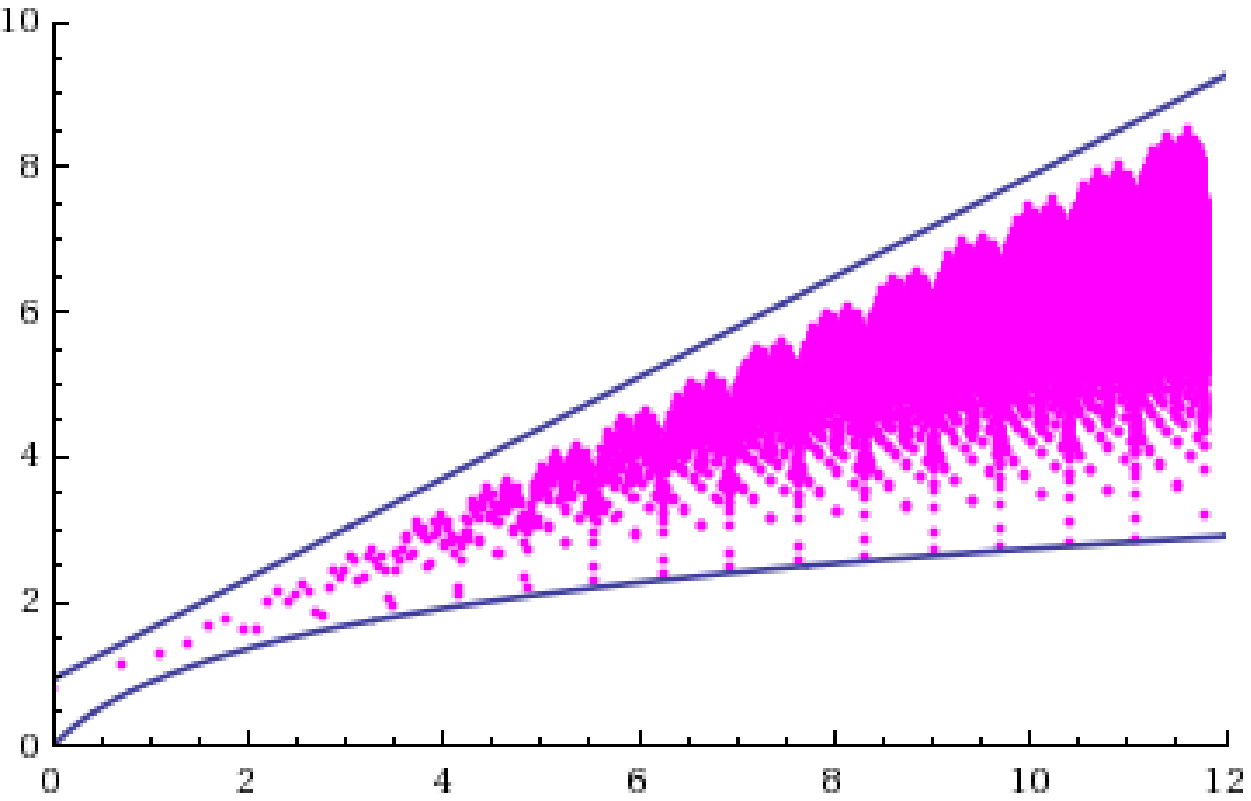}
        \put(7,55){\tiny $\log\|b_k\|$}
        \put(80,-5){\tiny $\log k$}
        \put(20,45){\tiny $y=\log_2\alpha_2 (x+1)+\log\frac{2\sqrt{2}}{5}$}
        \put(50,13){\tiny $y=\log(\frac{x}{\log2}+1)$}
      \end{overpic}
      &
      \begin{overpic}[tics=5,width=5.5cm]{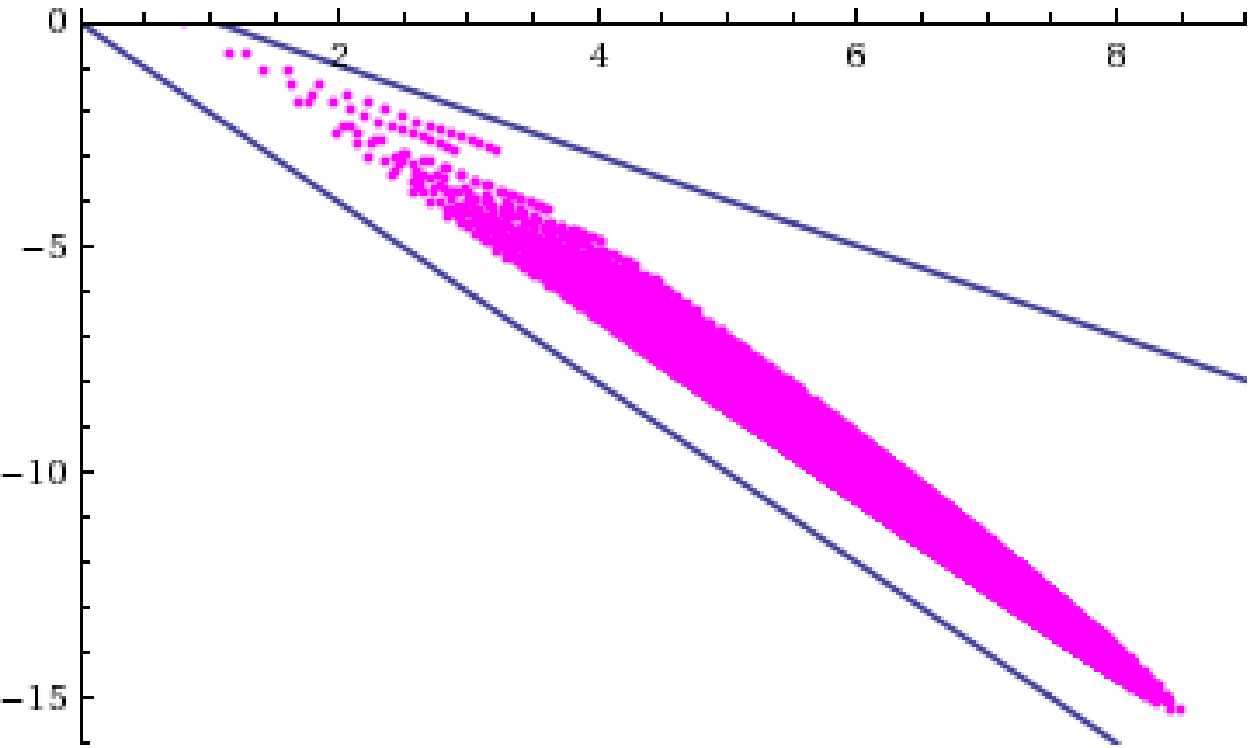}
        \put(10,42){\tiny $\log|S_I|$}
        \put(80,-5){\tiny $\log \|b_I\|$}
        \put(30,20){\tiny $y=-2x$}
        \put(70,45){\tiny $y=-x+1$}
      \end{overpic}
      \\
      (a)&(b)\\
    \end{tabular}
  \end{center}
  \caption{%
    \small 
    (a) log-log plot of the barycenters' norm vs. its index. The $b_I$ are arranged in the
    sequence naturally associated to the ordered tree $T(\cE)$, namely $b_I$ follows $b_J$
    if $|I|>|J|$ or, in case the multi-indices have the same order, the lowest index which 
    is different between $I$ and $J$ is bigger in $I$. Since there are $2^k$ nodes at the level
    $k$ the upper and lower bounds are evaluated using the fact that 
    $k+1\leq\|b_{2^k-1}\|\leq\frac{2\sqrt{2}}{\sqrt{5}}\alpha_2^{k+1}$ and therefore 
    $\log_2 k\leq\|b_k\|\leq\frac{2\sqrt{2}}{\sqrt{5}}(k+1)^{\log_2\alpha_2}$. (b) log-log plot of
    the length of the 1-simplices $S_I$ vs. the barycenters' norm. 
  }
  \label{fig:plots-n1}
\end{figure}
The growth rate of the Fibonacci sections here is given by the Golden Ratio $\alpha=(1+\sqrt{5})/2$ and
for the norm of the sections' barycenters we have the inequalities $k+2\leq\|b_k\|_\infty$ for the
slowest section and $\|b_k\|_\infty\leq \frac{2}{\sqrt{5}}\alpha^{k+1}$ for the fastest 
(see fig.~\ref{fig:plots-n1}(a)).
In particular the components of the two root Fibonacci sections are exactly the Fibonacci numbers:
e.g. taking $b_1=b(Z_\cE)=(1,2)$ and $b_2=b(Z_{\{e_1,e_1+e_2\}})=(2,3)$ we have that $b_3=(3,5)$, 
$b_4=(5,8)$ and so on.

Asymptotics of bodies have no meaning here but still we can say something about the asymptotics 
of the lengths of the 1-simplices constituting the binary tree $T(\cE)$. 
Indeed if $\cE'=\{ae_1+be_2,ce_1+de_2\}$ with $e_i=(x_i,y_i)$ then, in the chart $y=1$, 
$$
\vbox{\halign{\hfill$\displaystyle #$\hfill\cr
\mu(\cE')=d([ae_1+be_2],[ce_1+de_2])=\Big|\frac{ax_1+bx_2}{ay_1+by_2}-\frac{cx_1+dx_2}{cy_1+dy_2}\Big|\cr
=\frac{|ad-bc|\cdot|x_1y_2-x_2y_1|}{(ay_1+by_2)(cy_1+dy_2)}=\frac{1}{(ay_1+by_2)(cy_1+dy_2)}=\frac{1}{(a+b)(c+d)}\cr
}}
$$
where $|x_1y_2-x_2y_1|=1$ is the surface of the parallelogram corresponding to $\cE$ and 
$|ad-bc|=1$ because of the way the algorithm produces the new bases. 
In our concrete case $y_i=1$ and therefore
$$
\frac{1}{\|b\|^2_\infty}\leq \mu(\cE') =\frac{1}{\|b\|_\infty}(\frac{1}{a+b}+\frac{1}{c+d})\leq\frac{2}{\|b\|_\infty}
$$
Numerical illustrations of this pair of inequalities are shown in fig~\ref{fig:plots-n1}(b).
%
\subsection{The case n=2}
\label{sec:n2}
This is the only case where the polytopes corresponding to the bases and to the bodies
are of the same kind, namely triangles.
The algorithm that produces the fractal reduces here to the following:
\begin{alg}
\label{alg}
  \item{1.} On the three edges of the triangle $\Delta$ with vertices $\{[e_i]\}_{i=1,2,3}$ 
    select the three points $f_1=[e_2+e_3]$, $f_2=[e_3+e_1]$, $f_3=[e_1+e_2]$;
  \item{2.} subtract from $\Delta$ the interior of the triangle $Z$ (the ``body'' of $\Delta$)
    with vertices $\{f_1,f_2,f_3\}$;
  \item{3.} repeat recursively the algorithm on each of the three triangles that are left after
    the subtraction.
\end{alg}
Note that no two bodies have in common more than a point, i.e. they meet transversally, 
so the set $F(\cE)$ is never empty and actually it contains uncountably many points;
countably many of them can be explicitly evaluated through Fibonacci sections 
of the ternary tree $T(\cE)$.
\par
Consider for example the case 
$$\cE=\{e_1=(1,0,0),e_2=(0,1,0),e_3=(0,0,1)\}\,.$$
As shown in Theorem~\ref{thm:bar} 
the barycenter of every body triangle is the vector sum of the barycenters of the three body 
triangles it touches with its vertices (note that by construction no two bodies have a vertex 
in common) and when a body touches one of the sides of the root triangle $S(\cE)$ then we sum 
instead the vectors $(-1,1,1)$, $(1,-1,1)$ and $(1,1,-1)$ in correspondence respectively with
the sides opposed to the vertices $[e_1]$, $[e_2]$ and $[e_3]$.
The barycenters of one of the six root Fibonacci sections are determined by the first 
elements 
$$b_{-3}=(1,1,-1), b_{-2}=(1,-1,1),b_{-1}=(-1,1,1)$$
so that the generic element of the section is given by $b_k=(a_{k-2},a_{k-1},a_k)$, where
$a_k$ is the sequence of Tribonacci numbers with initial conditions $a_{-2}=1$, $a_{-1}=1$,
$a_0=1$. The expression of the generic term is given by
$$
a_k=
\frac{(1-\beta)(1-\bar\beta)}{(\alpha-\beta)(\alpha-\bar\beta)}\alpha^k+
\frac{(1-\alpha)(1-\bar\beta)}{(\beta-\alpha)(\beta-\bar\beta)}\beta^k+
\frac{(1-\alpha)(1-\beta)}{(\bar\beta-\beta)(\bar\beta-\alpha)}\bar\beta^k
$$
where $\alpha$, $\beta$ and $\bar\beta$ are the roots of the Tribonacci equation
$x^3=x^2+x+1$. Since $|\beta|<\alpha$ we have that $\|b_k\|_\infty\leq \frac{3(\alpha^2-1)}{3a^2-2a-1}\alpha^k$ 
and the limit point (see fig.~\ref{fig:spiral}) is $(1:\alpha:\alpha^2)$.
Note that all barycenters of this sequence lie on the ``projective Tribonacci spiral'' 
$$
\gamma(t)=[a(t-2):a(t-1):a(t)]
$$
where $a(t)$ is the trivial analytical extension of the $a_k$ sequence.
The fractal is invariant with respect to the projective transformations
$$
\vbox{\halign{\hfill$\displaystyle #$&$\displaystyle #$\hfill\cr
\psi_1([h_1:h_2:h_3])=&[h_1+h_2+h_3:h_2:h_3]\cr
\psi_2([h_1:h_2:h_3])=&[h_1:h_1+h_2+h_3:h_3]\cr
\psi_3([h_1:h_2:h_3])=&[h_1:h_2:h_1+h_2+h_3]\cr
}}
$$
so by applying any finite composition of them we obtain countably many explicit points of $F(\cE)$.

\begin{figure}
  \begin{center}
    \includegraphics[viewport=160 540 495 780,clip,width=15cm]{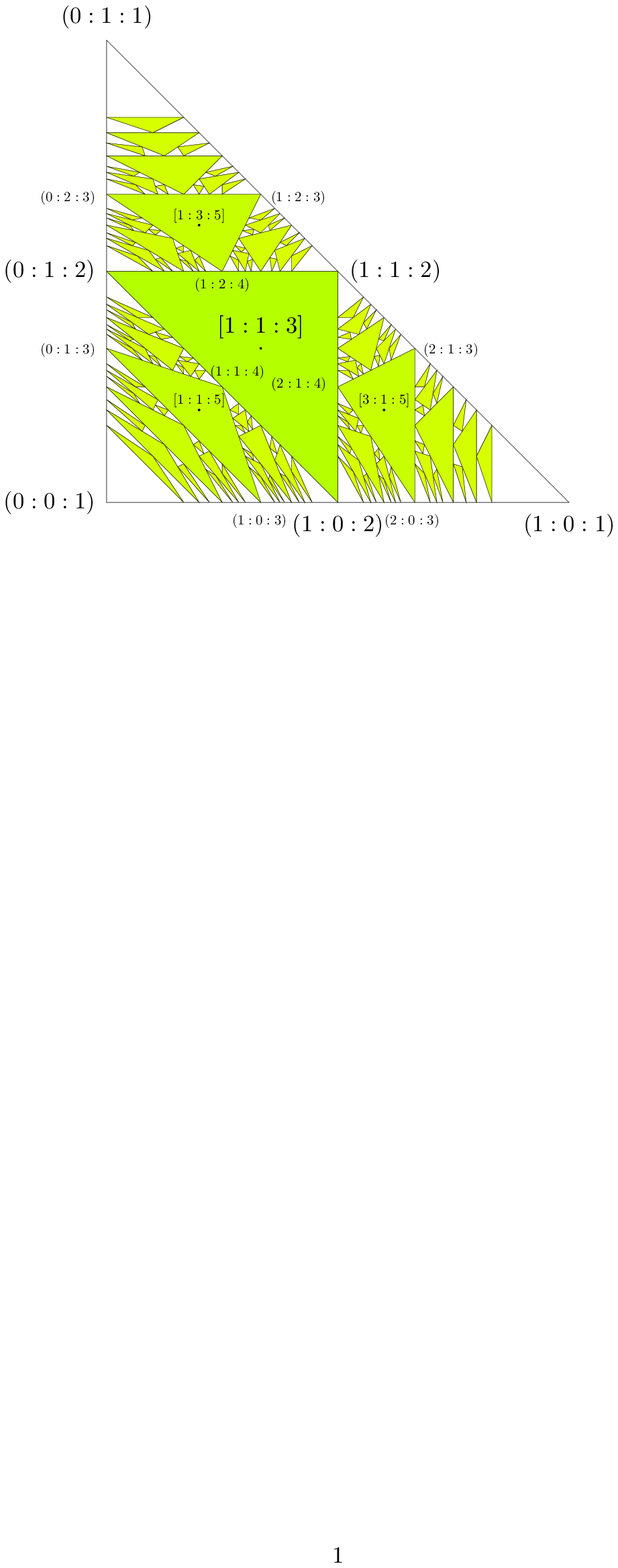}
    \caption{%
      \small
      Plot of $F^5(\cE)$, namely of the bodies up to the forth recursion level, for 
      $\cE=\{(1,0,1),(0,1,1),(0,0,1)\}$ in the $h^3=1$ projective chart of $\bR \hbox{P}^2$.
      Bodies are colored in green, so the points of $F^5$ are the white ones. 
      The homogeneous coordinates of the vertices of 
      the first and second level bodies are shown together with the body's barycenters, for which
      we used the square brakets for sake of clarity. Note that barycenters can be obtained in three 
      ways: (i) by summing the barycenters of the three bodies touched by the vertices -- note that
      in case a vertex touches a root simplex edge then the following should be used: $(1,1,1)$ 
      for the edge opposite to $[0:0:1]$, $(1,-1,1)$ for the one opposite to $[1:0:1]$ and $(-1,1,1)$ 
      for the one opposite to $[0:1:1]$; (ii) by summing the coordinates of the vertices of the triangle 
      that generated the body; (iii) by summing the coordinates of the vertices of the body and
      dividing them by 2 -- this corresponds to the fact the volume associated to the 
      basis corresponding to the vectors $\{e_1+e_2,e_2+e_3,e_3+e_1\}$ is double with respect
      to the basis $\{e_1,e_2,e_3\}$.
    }
    \label{fig:fr-n2}
  \end{center}
\end{figure}
The slowest sections in the barycenters' norms growth is, modulo indices permutations,
$$t_k=\cE_{\underbrace{1,\cdots,1}_{k}}=\{e^{(k)}_1=e_1,e^{(k)}_2=e_2+ke_1,e^{(k)}_3=e_3+ke_1\}$$
for which $b_k=(2k+3,1,1)$ and therefore $\|b_k\|_\infty=2k+3$.

\begin{figure}
  \label{fig:asymp-n2}
  \begin{center}
    \begin{tabular}{cc}
      \begin{overpic}[tics=5,width=5.5cm]{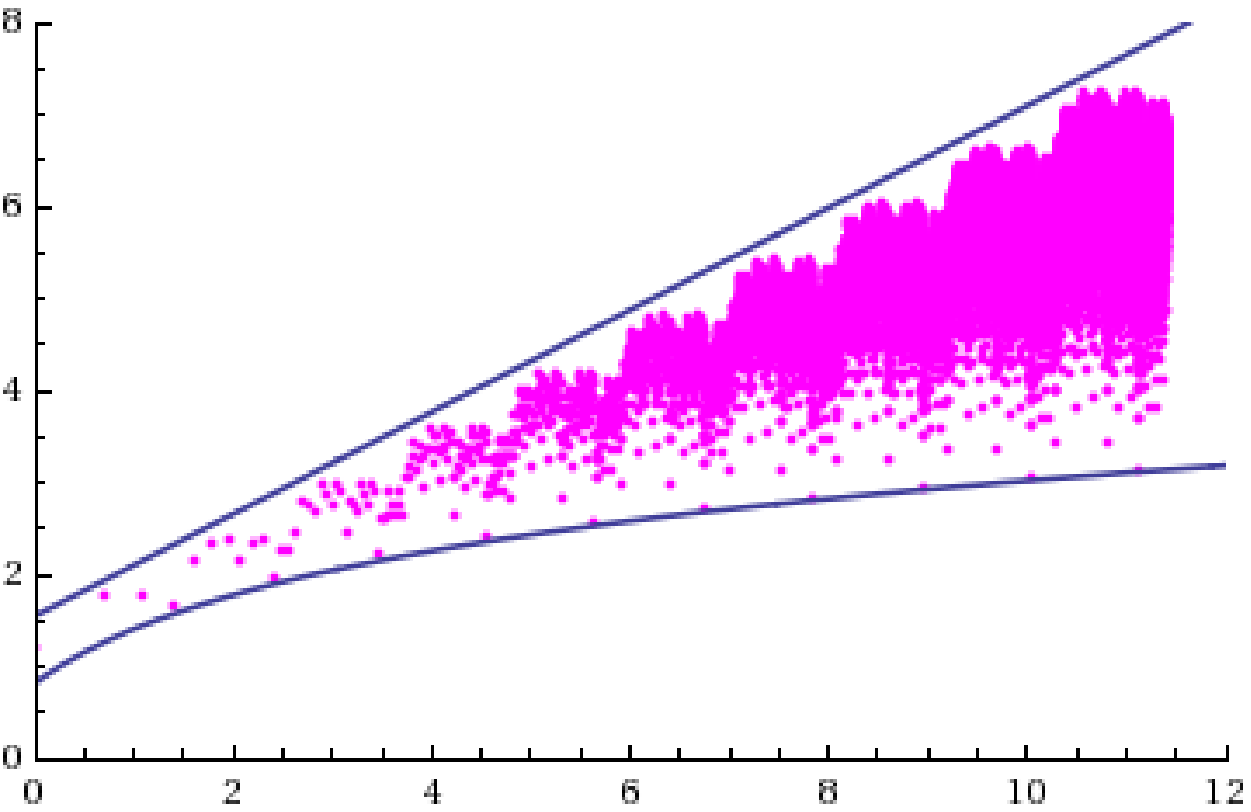}
        \put(6,55){\tiny $\log\|b_k\|$}
        \put(80,-3){\tiny $\log k$}
        \put(17,45){\tiny $y=\log_3\alpha_3 x+\log A$}
        \put(50,17){\tiny $y=\log(\frac{2x}{\log3}+\log_312)$}
      \end{overpic}
      &
      \begin{overpic}[tics=5,width=5.5cm]{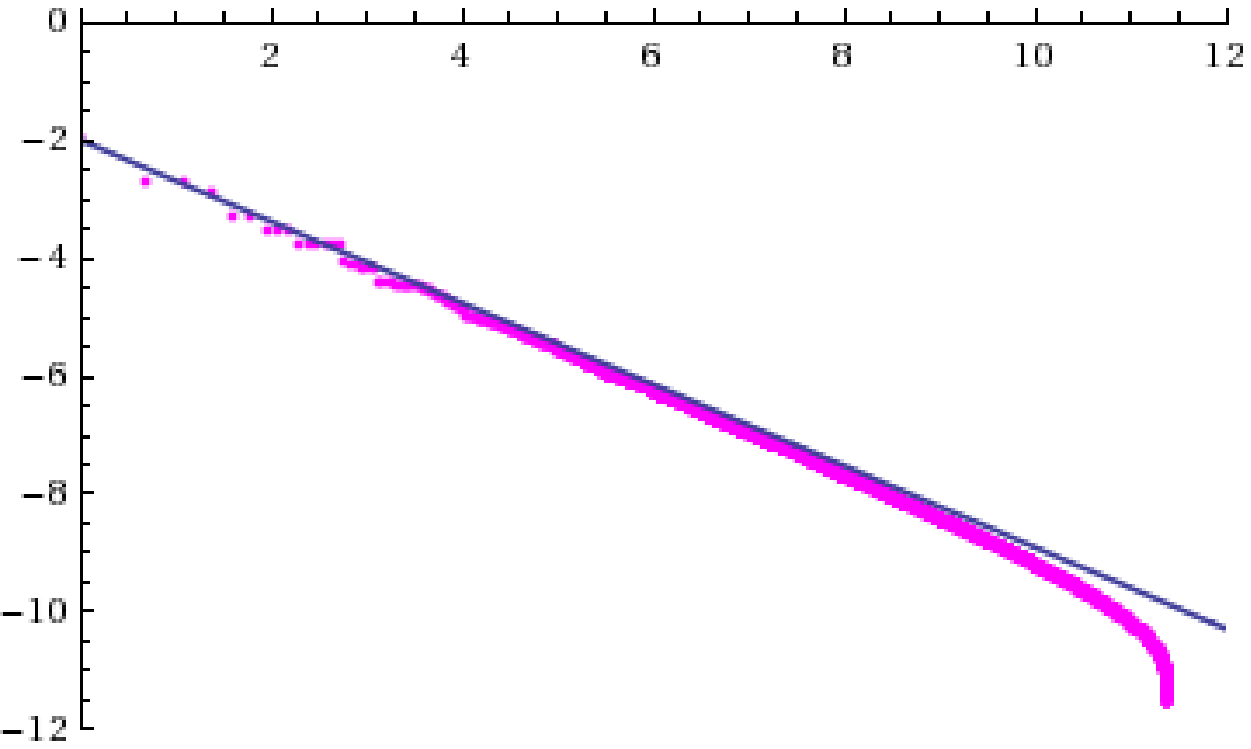}
        \put(9,53){\tiny $\log\|\rho_k\|$}
        \put(80,-3){\tiny $\log k$}
        \put(50,36){\tiny $y=-.69k-2$}
      \end{overpic}
      \\
      (a)&(d)\\
      \begin{overpic}[tics=5,width=5.5cm]{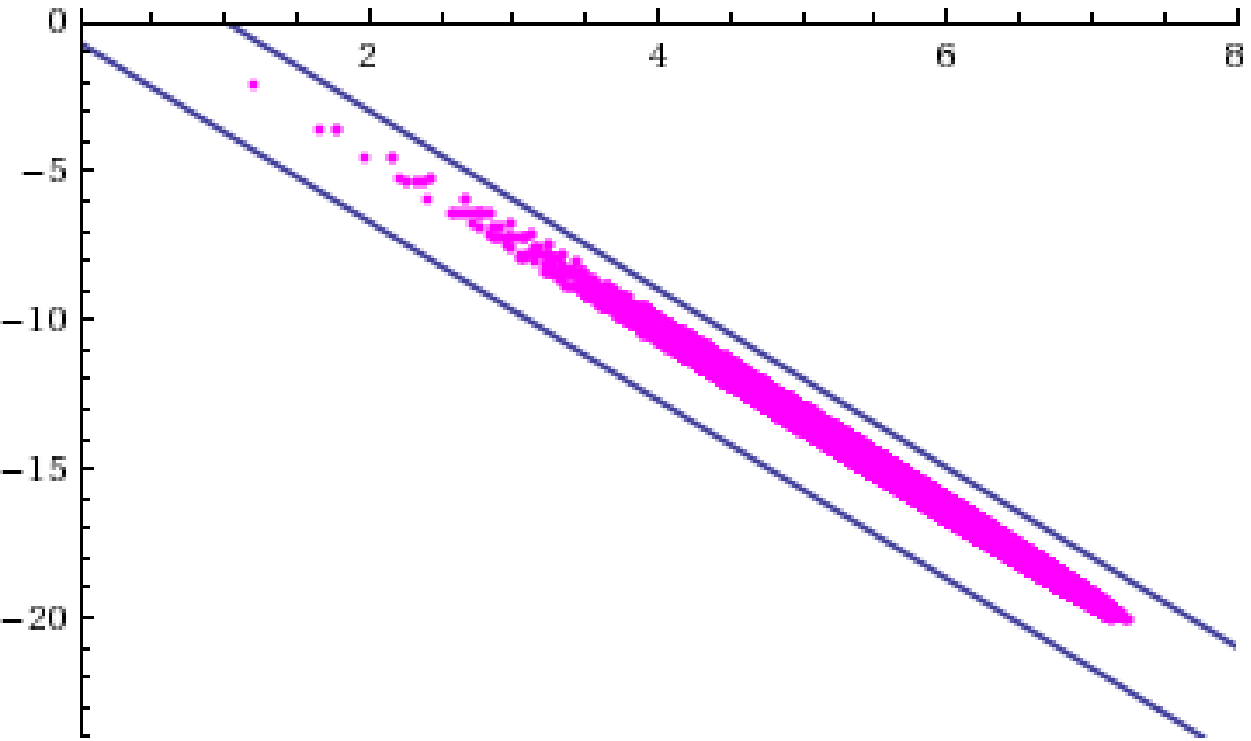}
        \put(9,37){\tiny $\log\mu(Z_I)$}
        \put(80,-3){\tiny $\log \|b_I\|$}
        \put(50,42){\tiny $y=-3x+\log(8\sqrt{3})$}
        \put(34,16){\tiny $y=-3x-\log2$}
      \end{overpic}
      &
      \begin{overpic}[tics=5,width=5.5cm]{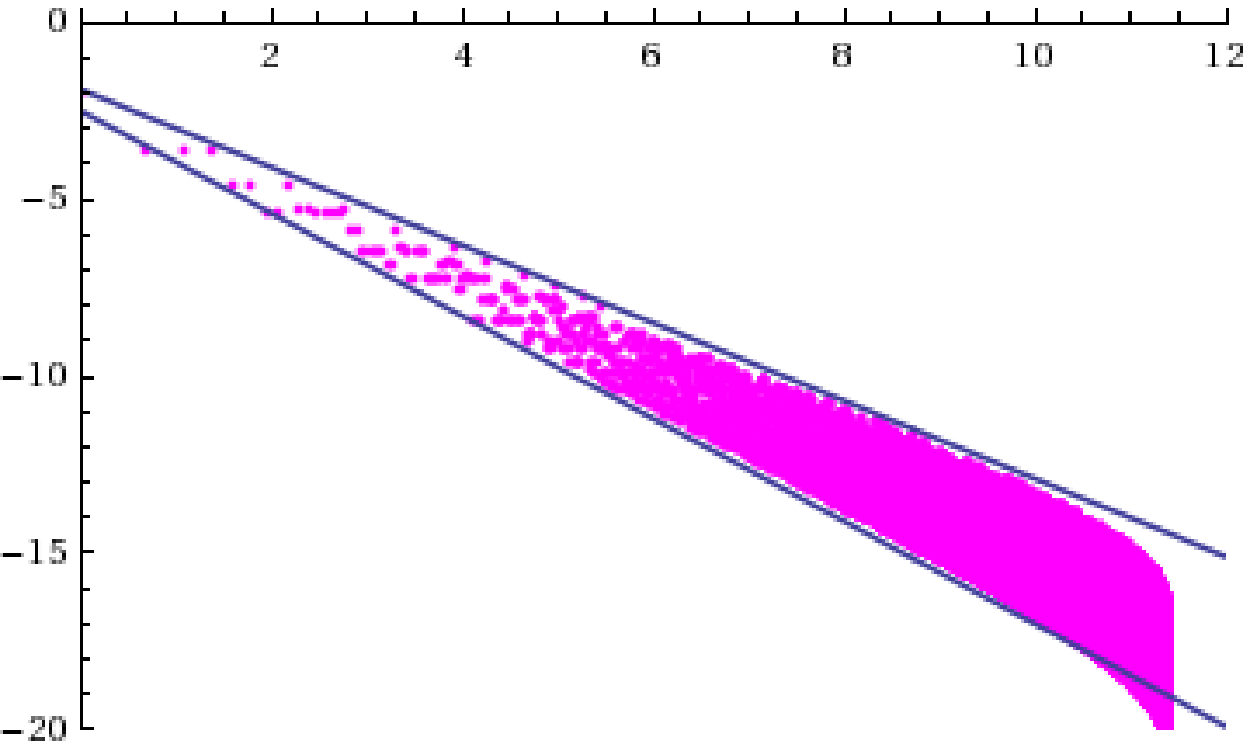}
        \put(9,37){\tiny $\log\mu(Z_k)$}
        \put(80,-3){\tiny $\log k$}
        \put(50,42){\tiny $y=-1.1x-1.9$}
        \put(32,16){\tiny $y=-1.45x-2.5$}
      \end{overpic}
      \\
      (b)&(e)\\
      \begin{overpic}[tics=5,width=5.5cm]{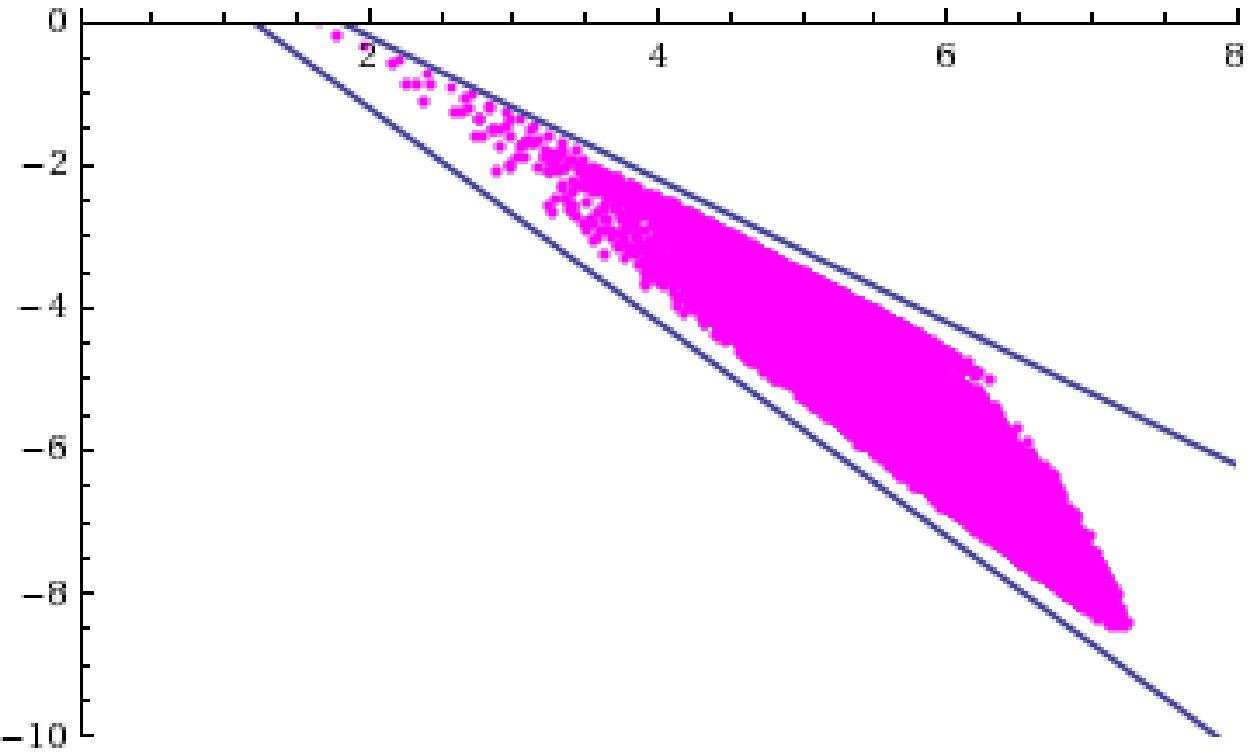}
        \put(9,37){\tiny $\log|Z_I|$}
        \put(80,-3){\tiny $\log \|b_I\|$}
        \put(70,42){\tiny $y=-x+\log6$}
        \put(36,18){\tiny $y=-\frac{3}{2}x+1.8$}
      \end{overpic}
      &
      \begin{overpic}[tics=5,width=5.5cm]{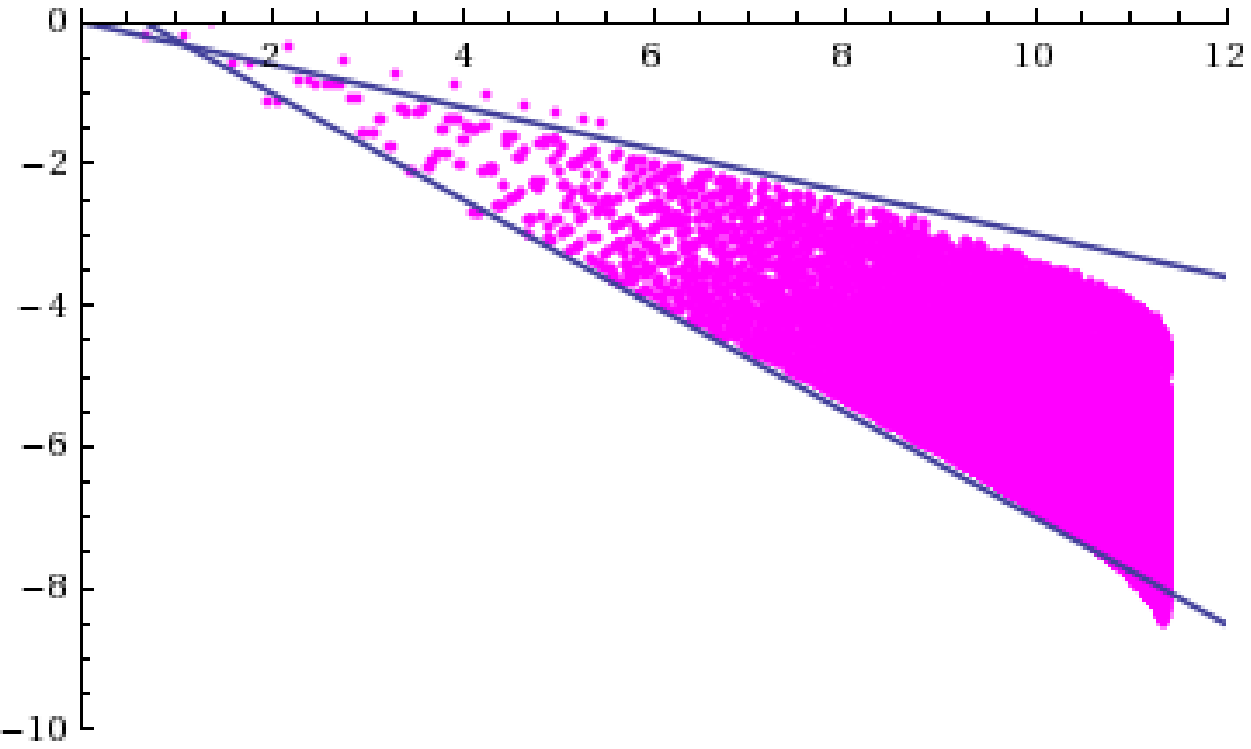}
        \put(9,37){\tiny $\log|Z_k|$}
        \put(80,-3){\tiny $\log k$}
        \put(68,48){\tiny $y=-.3x$}
        \put(31,28){\tiny $y=-.75x+.5$}
      \end{overpic}
      \\
      (c)&(f)\\
    \end{tabular}
  \end{center}
  \caption{%
    \small
    Log-log plots for the main quantities in the $n=2$ case for $\cE=\{(1,0,1),(0,1,1),(0,0,1)\}$ 
    (a) Barycenters norms vs indices -- as explained in fig.~\ref{fig:plots-n1} the $b_k$ are ordered 
    according to the natural order induced by the tree, so that 
    $2k+3\leq\|b_\frac{3^k-1}{2}\|\leq\sqrt{3}\alpha_3^3\frac{(1-\beta_3)(1-\bar\beta_3)}{(\alpha-\beta_3)(\alpha-\bar\beta_3)}\alpha_3^k=A\alpha_3^k$
    and therefore
    $\frac{2}{\log3}\log k+\log_312 \leq \|b_k\| \leq Ak^{\log_3\alpha_3}$.
    (b) Bodies' volumes vs barycenters norms and (c) bodies' diameters vs barycenters norms -- the
    lines bounding the numerical data come immediately from the inequalities in Theorem~\ref{thm:bodyVol-n2}.
    For the next three plots no exact formulae are known so the lines shown represent just an
    interpolation of the numerical data. 
    (d) Radii of the circles inscribed in the bodies vs k after sorting the radii in descending order.
    (e) Areas of the bodies and (f) their diameter sorted according with their radii.
  }
\end{figure}

In figs.~\ref{fig:asymp-n2}(a-c) we show the numerical results for the asymptotic
behaviour of the barycentric norms and the bodies' surfaces and diameters.

Note that in this particular case Conjecture~\ref{conj:Diam} is known to be true
through an indirect proof. Indeed, this fractal comes up naturally in the study of 
the asymptotics of plane sections of periodic surfaces, which in turn comes from the 
problem of the motion of quasi-electrons under a strong magnetic field (see~\cite{NM03}
for a detailed account), in the particular case of the regular triply-periodic skew 
polyhedron $\{4,6|4\}$~\cite{DD09}.
In that setting
the basis is 
$$\cE_\cC=\{e_1=(1,0,1),e_2=(0,1,1),e_3=(1,1,0)\}$$
and the barycenter $b$ of a body $Z$ represents a homological discrete ``first integral'' 
of a Poisson dynamical system which dictates the asymptotic directions of the plane sections
in the following way: the open sections obtained by cutting the polyhedron with planes
perpendicular to every direction $\omega\in Z$ are all strongly asymptotic to the direction
``$\omega\times b$''.
It is a general theorem of that theory the fact that the diameter of a body $Z$ 
is bounded by $C/\|b(Z)\|$ where $C$ is a constant depending only on the 
surface~\cite{DeL05a}, which then establishes the following theorem for this $n=2$ case:
\begin{theorem}
  \label{thm:bodyVol-n2}
  Let $Z$ be a body in $T(\cE_\cC)$ with area $\mu(Z)$ (where $\mu$ is the same measure
  used in Theorem~\ref{thm:bodyVol}), diameter $|Z|$ and barycenter $b$. 
  Then 
  the following inequalities hold asymptotically:  
  $$
  \frac{1}{^4\sqrt{3}\|b\|^{\frac{3}{2}}}\leq |Z| \leq\frac{6}{\|b\|},\frac{1}{2\|b\|^3}\leq \mu(Z) \leq\frac{12\sqrt{3}}{\|b\|^3}\,.
  $$ 
\end{theorem}
\begin{proof}
  The inequality for the area of $Z$ is just the restriction of Theorem~\ref{thm:bodyVol} to $n=2$ together
  with the fact that $\|b\|_\infty\leq\|b\|\leq\sqrt{3}\|b\|_\infty$. The right hand side for the diameter comes 
  from the general theory of plane sections of a triply periodic surface that, applied to
  this particular case, states~\cite{DeL05a} that the distance between the barycenter and the bodies' vertices 
  is bounded by $3/\|b\|$, where the 3 is the double of the area of the basic cell of the periodic
  surface cited above in this section. The left hand side comes simply from the fact that a triangle of area $a$ 
  cannot have a diameter smaller than $\sqrt{2a/\sqrt{3}}$.
\end{proof}
Being unable to evaluate analytical bounds for the Hausdorff dimension $d_\cC$ of $F(\cE_\cC)$, we compute 
numerically four different quantities that may give hints on whether $d_\cC$ is integer or not (the
non-integrality of $d_\cC$ would confirm a general conjecture by Novikov~\cite{NM03}).


%
First we get a direct upper bound for the Hausdorff dimension by counting the smallest number 
of squares of side $\epsilon=2^{-l}$, $l=0,\cdots,12$, needed to cover $F^{12}$, i.e. the union of
all bodies up to the 12-th order of recursion; as shown in fig.~\ref{fig:asymp-n2}(h), 
we get $d_\cC\lesssim1.7$.

Then we evaluate the {\sl Minkowsky dimension}, namely the limit
$$
2-\lim_{\epsilon\to0}\frac{\log V(F_\epsilon)}{\log \epsilon}
$$
where $V(F_\epsilon)$ is the surface of the $\epsilon$ neighborhood of $F$,
using the formula~\cite{Fal97,Gai06}
$$
V(F_\epsilon) = p\epsilon + \epsilon\sum_{i=1}^{k_\epsilon} p_i + A - \sum_{i=1}^{k_\epsilon} a_i + 
\epsilon^2(\pi-\sum_{i=1}^{k_\epsilon}\frac{p_i^2}{4a_i})
$$
where $k_\epsilon$ is the integer such that $\rho_{k_\epsilon+1}\leq\epsilon\leq\rho_{k_\epsilon}$,
$\rho_k$ is the radius of the inscribed circle to the body $Z_k$ and the bodies are
sorted in descending order with respect to the radii.
In fig.~\ref{fig:asymp-n2}(g) we show the numerical results we got by evaluating the volume 
of the neighborhoods of $\cE$ of radii $r_n=1.2^{-n}$ for $n=1,\cdots,50$, which suggests 
a Minkowsky dimension between $1.7$ and $1.8$.

Next, we evaluate numerically the growth rate of the radii after sorting them in
decreasing order (fig.~\ref{fig:asymp-n2}(d)) and then the corresponding bounds for the
bodies areas (fig.~\ref{fig:asymp-n2}(e)) and diameters (fig.~\ref{fig:asymp-n2}(f)).
In this case we obtain that $\epsilon\asymp k^{-.69}$, $Ak^{-1.45}\leq a_k \leq Bk^{-1.1}$
and $A'k^{-.75}\leq p_k \leq B'k^{-.3}$, so that
$A''\epsilon^{.65}\leq V_\epsilon \leq B''\epsilon^{.145}$. From this we get a second
evaluation, compatible but much looser, for the Minkowsky dimension: $1.35\leq \dim_MF(\cE_\cC)\leq1.86$.

Finally, we use Corollary~\ref{cor:bardensity} and evaluate the box-counting dimension
of the set $B=\{[b_I(\cE_\cC)]\}\subset\RPt$, namely the set of barycenters of the bases
in the tree $T(\cE_\cC)$ considered as points in the projective plane. Since the closure 
of $B$ is the union of $F(\cE_\cC)$ with the (one-dimensional) boundaries of the bodies, 
a dimension higher than one must be due to the points in $F(\cE_\cC)$. As shown in
Fig.~\ref{fig:labelsDim}, the  dimension appear to be about 1.69. Notice 
that applying this method to the Sierpinki triangle, whose dimension is $d_S=\log3/\log2$,
gives the correct approximation to the third digit $d_S\simeq1.59$.

In conclusion, the four evaluations are in excellent agreement with each other and 
indicate a non-integer Hausdorff dimension for this fractal, probably about 1.7; 
finding exact bounds 
would be nicer though since this would represent the first analytical confirmation 
of a conjecture of Novikov about the non-integer dimension of fractals coming from 
the theory of asymptotics of plane sections of triply-periodic surfaces. 
\subsection{The case n=3}
\label{sec:n3}
\begin{figure}
  \label{fig:n3}
  \begin{center}
    \begin{tabular}{cc}
      \includegraphics[width=5.5cm]{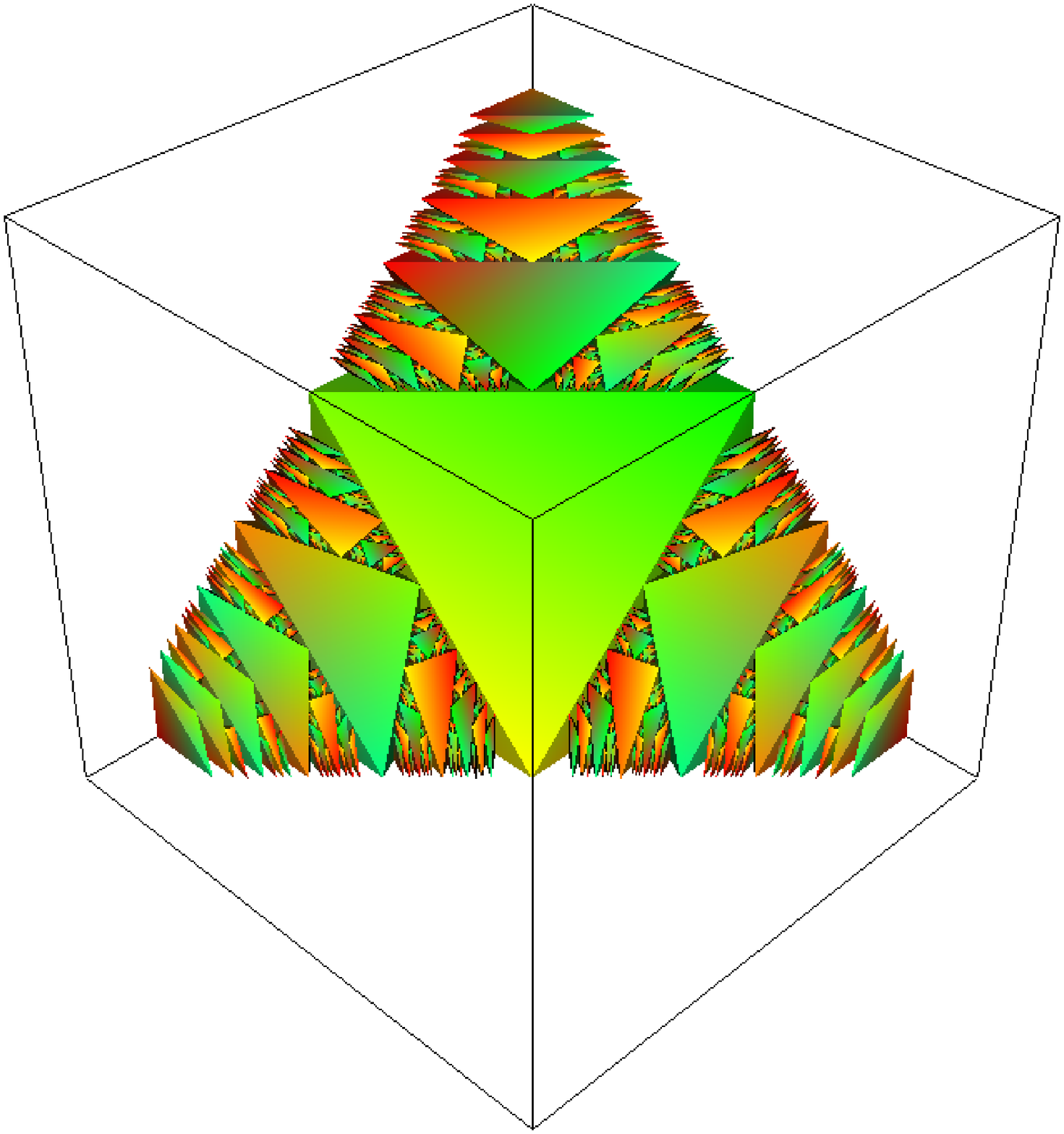}&
      \includegraphics[width=5.5cm]{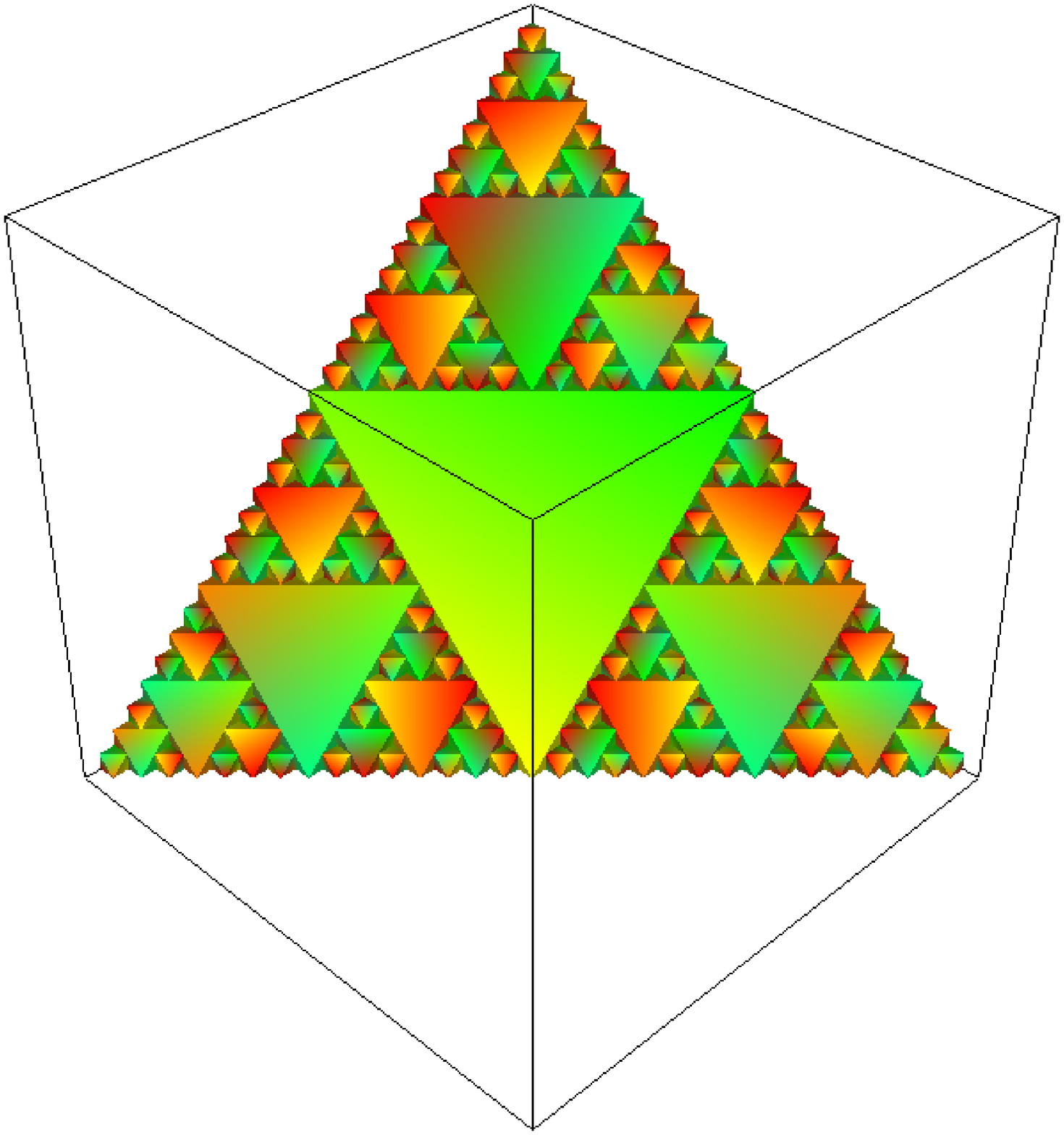}\\
      (a)&(b)\\
      \includegraphics[width=5.5cm]{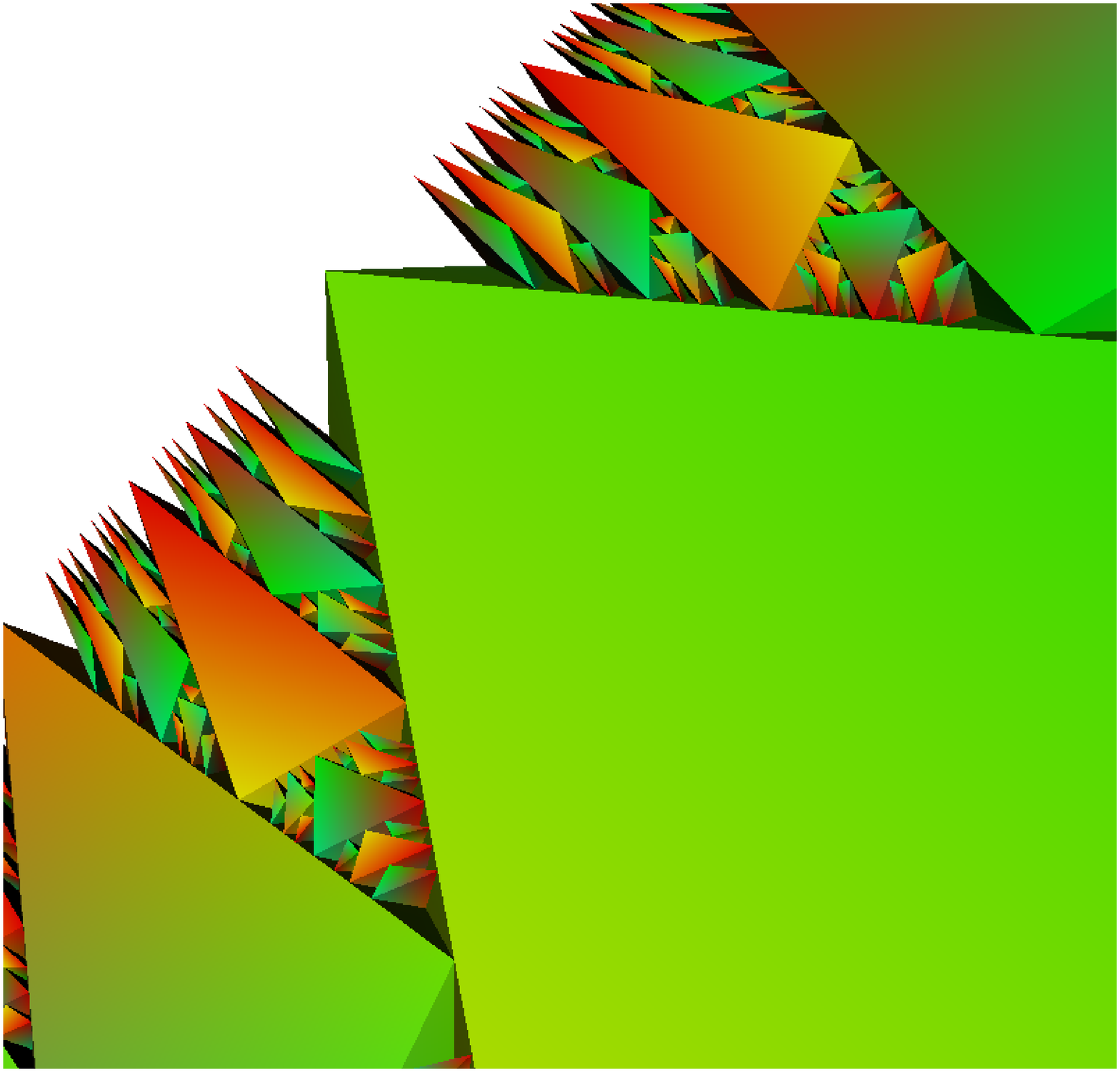}&
      \includegraphics[width=5.5cm]{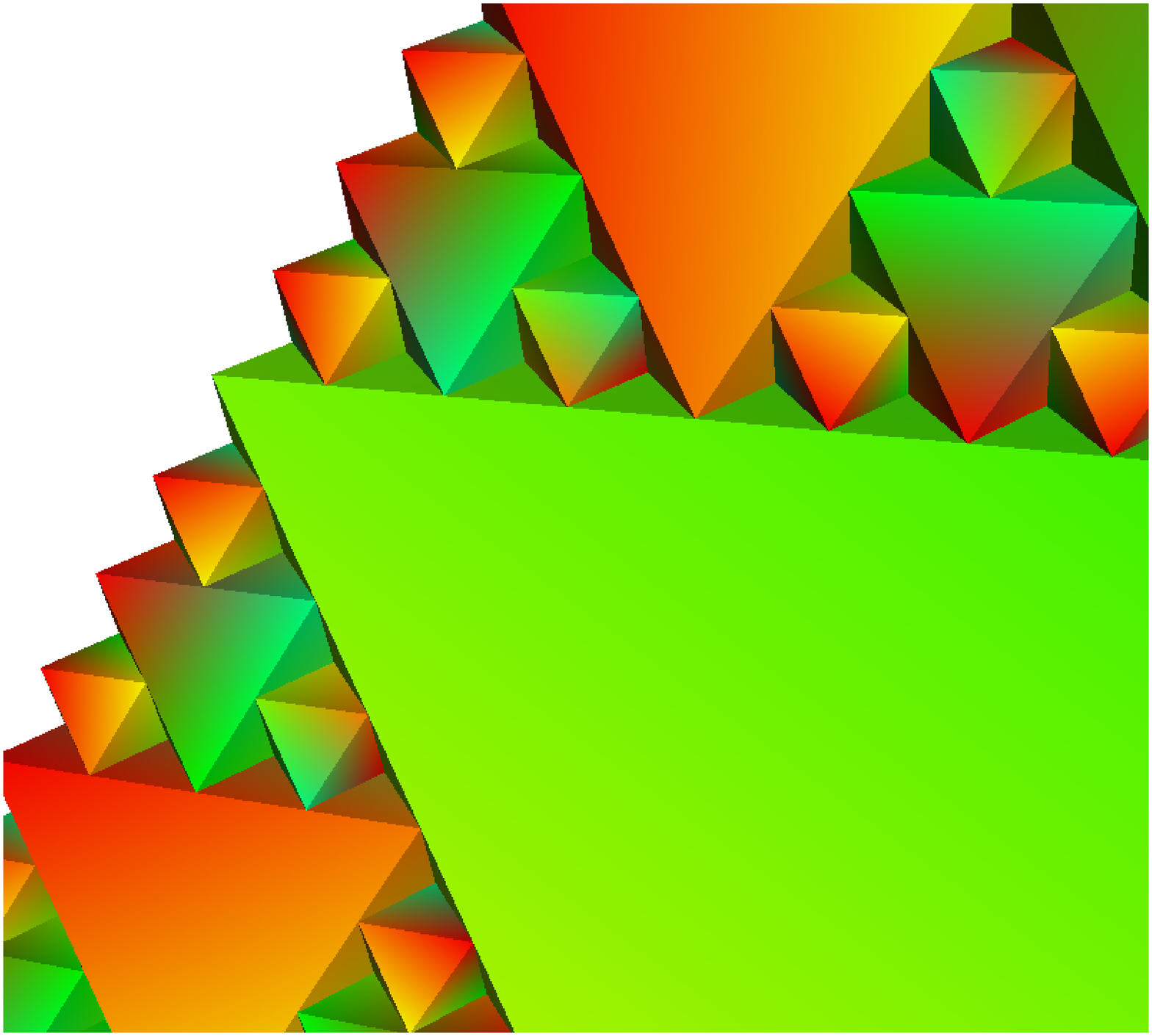}\\
      (c)&(d)\\
    \end{tabular}
  \end{center}
  \caption{%
    \small
    (a,c) Total view and detail of $F^6(\cE)$, namely of the bodies up to the fifth recursion level,
    for $\cE=\{(1,0,0,1),(0,1,0,1),(0,0,1,1),(0,0,0,1)\}$ 
    in the $h^4=1$ projective chart of $\bR \hbox{P}^3$. The bodies of $S(\cE)$ are shown, 
    up to the fifth recursion level, in red-green colors; the points of $F^6(\cE)$ are their 
    complement in the tetrahedron of vertices (in the chart $h^4=0$) $(1,0,0)$, $(0,1,0)$ and $(0,0,1)$.
    (b,d) Total view and detail, up to the forth recursion level, of the Sierpi\'nski tetrahedron.
  }
\end{figure}

When $n=3$, every body has 6 vertices: one for each edge of the tetradedron they belong to, 
and eight triangular faces, one for each face and one for each corner of the tetrahedron. 
Bodies that touch each other share a whole triangle (rather than a single point as in the 
$n=2$ case) in the following way: bodies can meet only on the faces that do not come from 
the tetrahedra $S_I$ and, on those faces, these shared triangles form a fractal of the $n=2$ 
kind (see Fig.~\ref{fig:n3}).

In the particular case of $\cE_T=\{(1,0,0,1),(0,1,0,1),(0,0,1,1),(0,0,0,1)\}$ the barycenter of
the root tetrahedron is $(1,1,1,4)$ and the volume inequalities translate in 
$$
\frac{2}{3\|b_I\|^4_\infty}\leq\mu(Z_I(\cE_T))\leq\frac{2\cdot3^3}{\|b_I\|^4_\infty}\,.
$$
In case of barycenters' norms we have $4+3k\leq\|b_k\|$ for the slowest tree section and 
$\|b_k\|\leq2\frac{(1-\beta)(1-\bar\beta)(1-\gamma)}{(\alpha_4-\beta)(\alpha_4-\bar\beta)(\alpha_4-\gamma)}\alpha_4^{k+4}$
for the fastest.

Numerical evaluations of the Hausdorff dimension are more cumbersome for $n=3$ because
the number of bodies grows very large after few iterations of the generating algorithm
(getting rather heavy on both CPU and RAM consumption) and their geometry gets much 
more complicated. 

\begin{figure}
  \label{fig:frdim-n3}
  \begin{center}
    \begin{tabular}{cc}
      \begin{overpic}[tics=5,width=5.5cm]{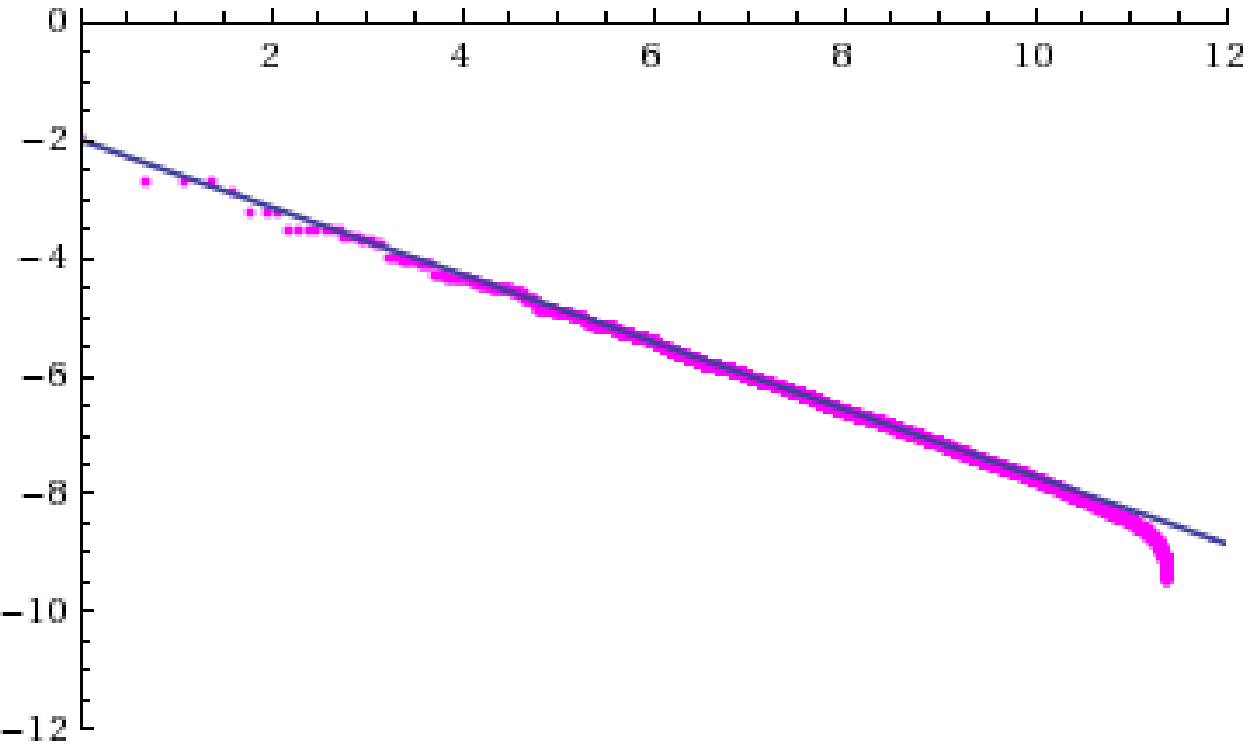}
        \put(9,53){\tiny $\log\|\rho_k\|$}
        \put(80,-3){\tiny $\log k$}
        \put(50,36){\tiny $y=-.57k-2$}
      \end{overpic}
      &
      \begin{overpic}[tics=5,width=5.5cm]{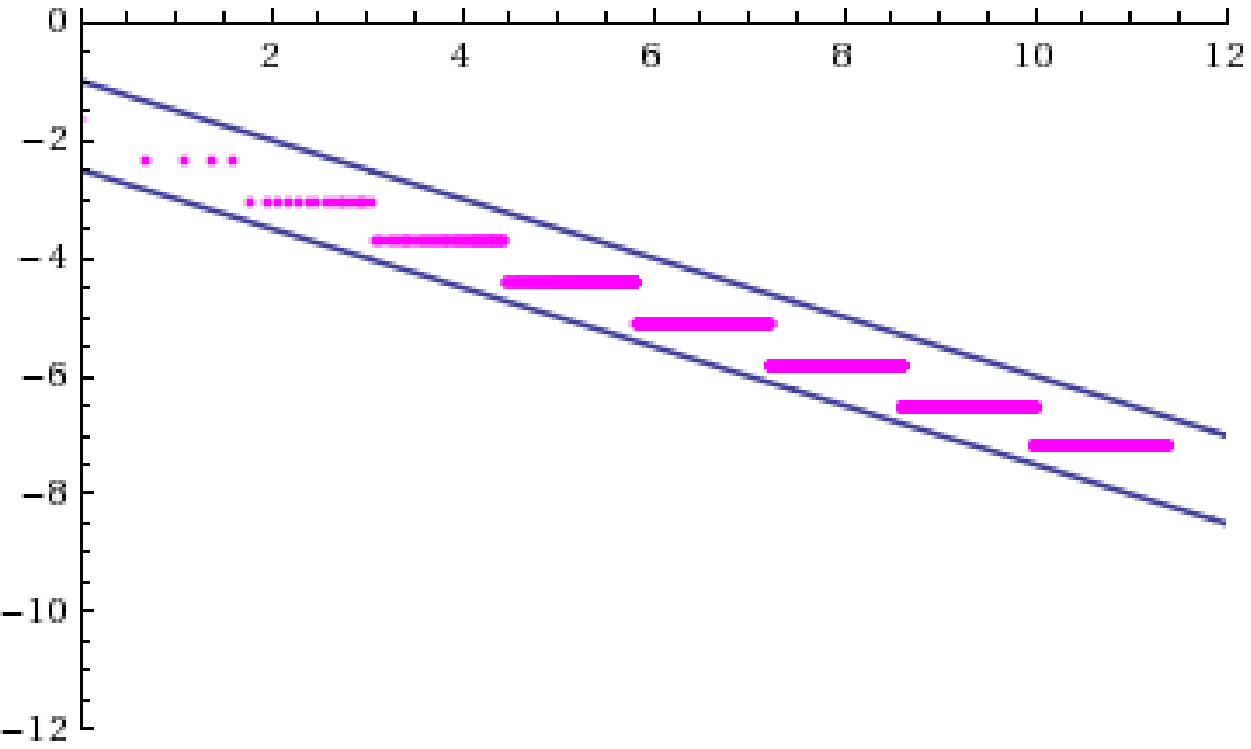}
        \put(9,53){\tiny $\log\|\rho_k\|$}
        \put(80,-3){\tiny $\log k$}
        \put(55,42){\tiny $y=-.5k-2.5$}
        \put(35,26){\tiny $y=-.5k-1$}
      \end{overpic}
      \\
      (a)&(d)\\
      \begin{overpic}[tics=5,width=5.5cm]{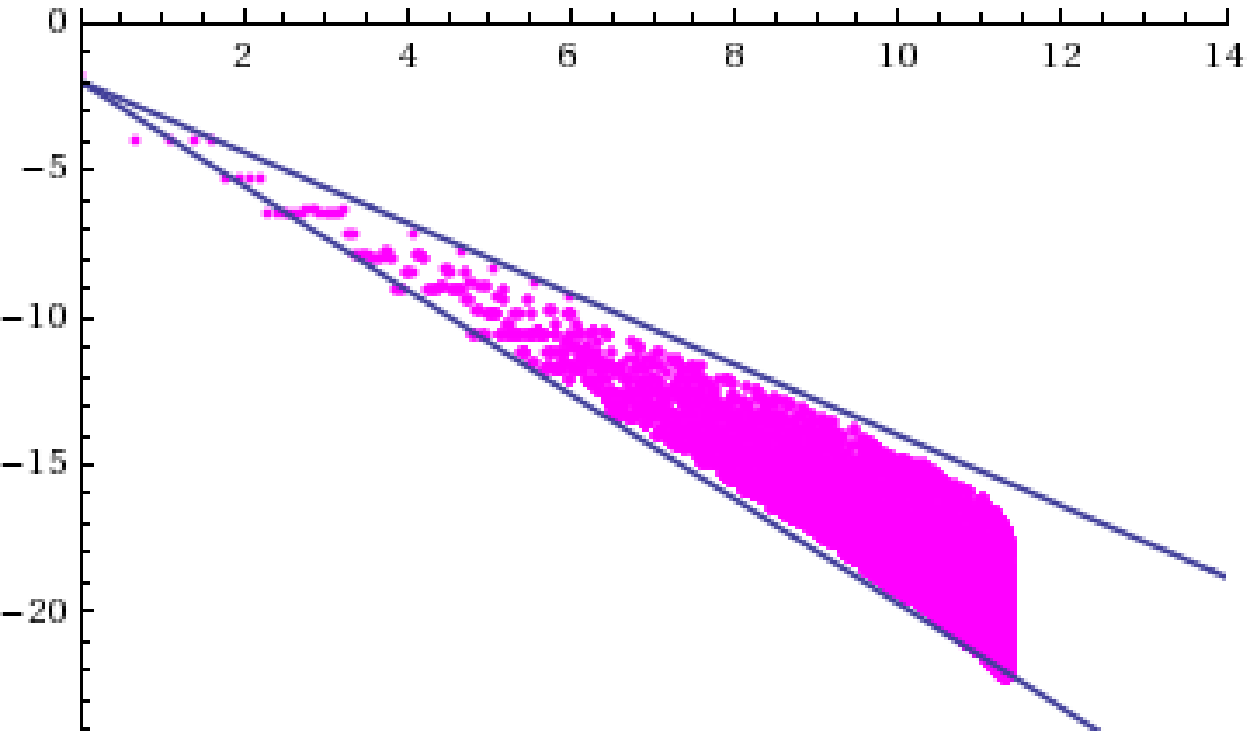}
        \put(9,37){\tiny $\log\mu(Z_I)$}
        \put(80,-3){\tiny $\log \|b_I\|$}
        \put(50,35){\tiny $y=-1.2x-2$}
        \put(28,16){\tiny $y=-1.77x-2$}
      \end{overpic}
      &
      \begin{overpic}[tics=5,width=5.5cm]{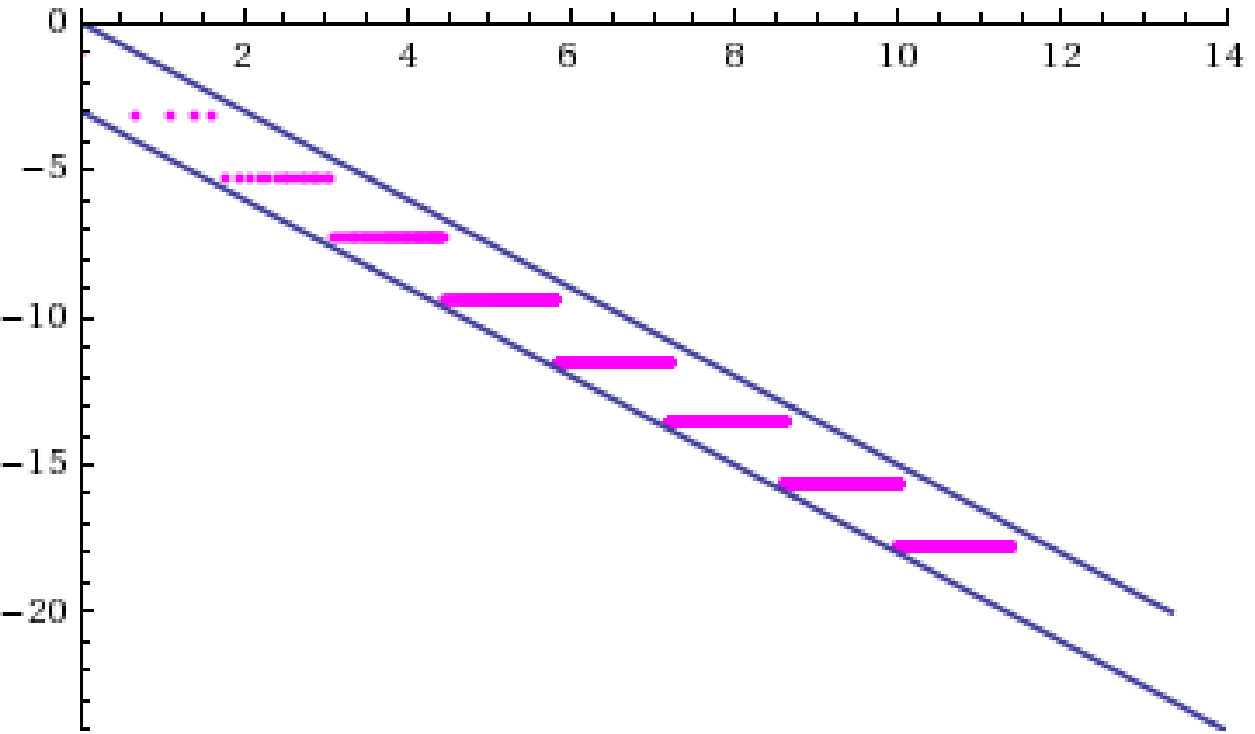}
        \put(9,37){\tiny $\log\mu(Z_I)$}
        \put(80,-3){\tiny $\log \|b_I\|$}
        \put(50,35){\tiny $y=-1.5x$}
        \put(28,16){\tiny $y=-1.5x-3$}
      \end{overpic}
      \\
      (b)&(e)\\
      \begin{overpic}[tics=5,width=5.5cm]{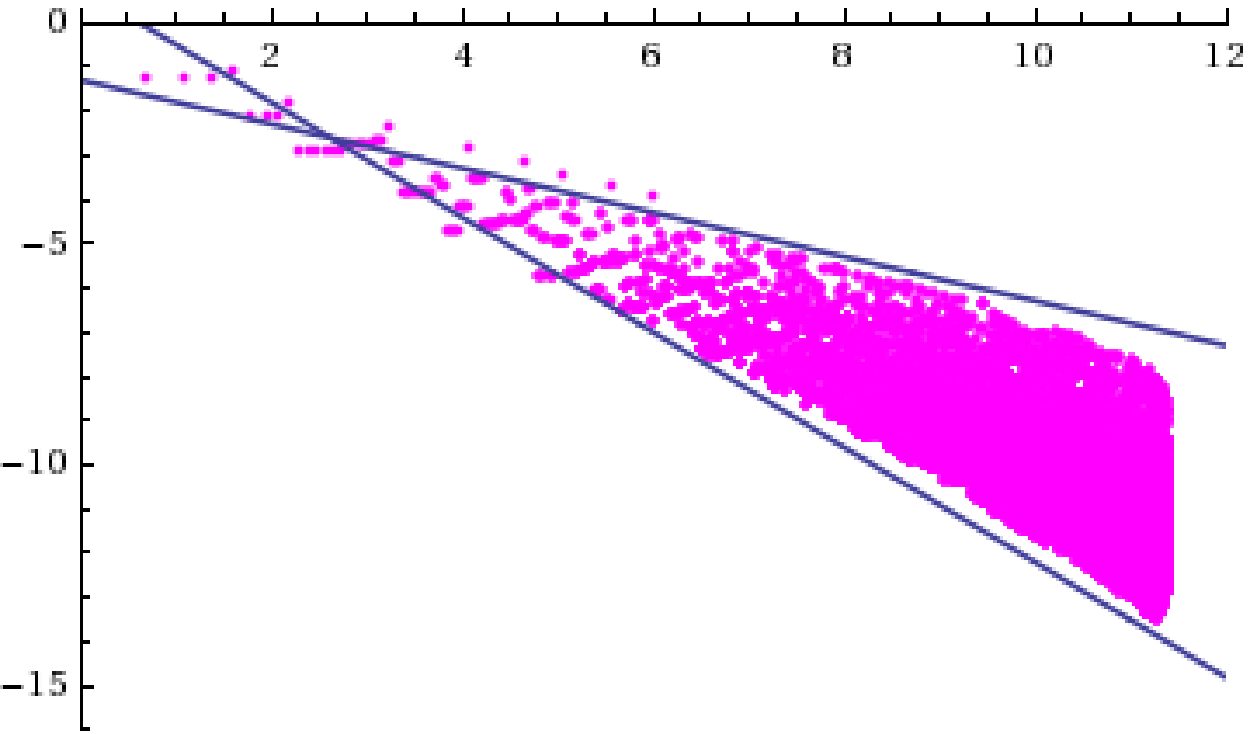}
        \put(9,37){\tiny $\log\mu(Z_I)$}
        \put(80,-3){\tiny $\log \|b_I\|$}
        \put(60,42){\tiny $y=-.5x-1.3$}
        \put(31,19){\tiny $y=-1.3x+.8$}
      \end{overpic}
      &
      \begin{overpic}[tics=5,width=5.5cm]{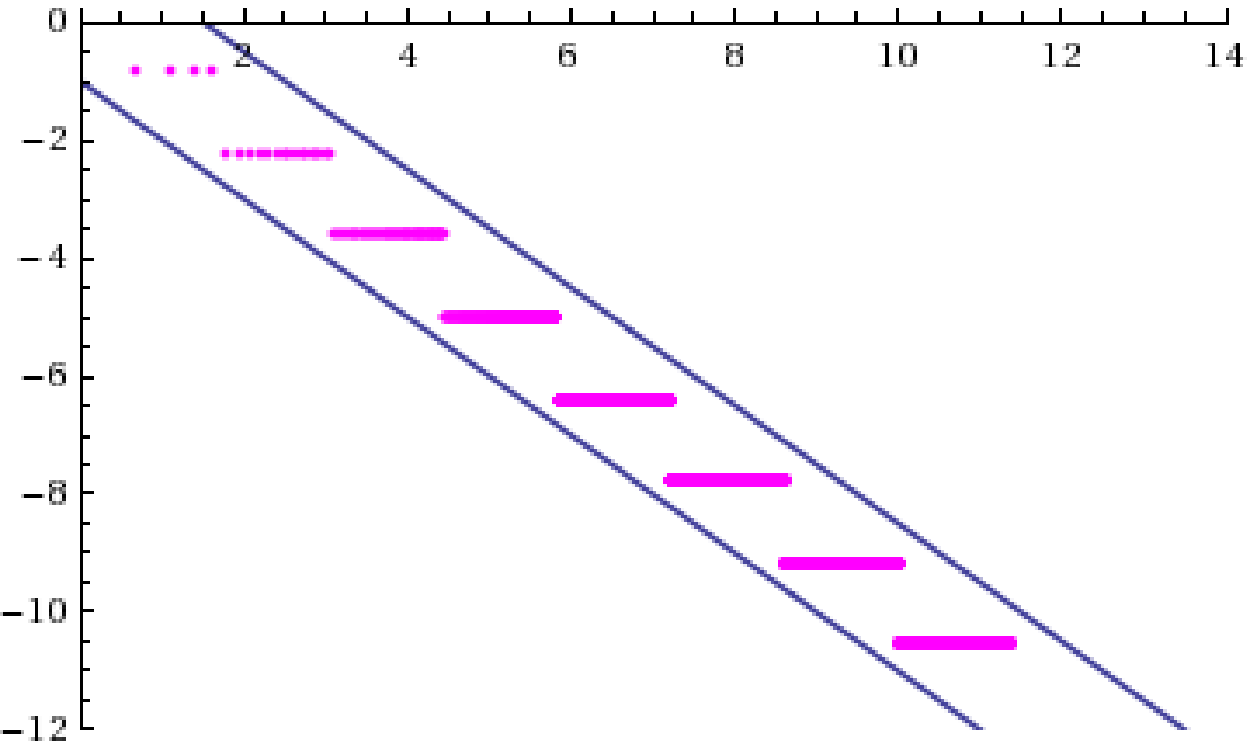}
        \put(9,37){\tiny $\log\mu(Z_I)$}
        \put(80,-3){\tiny $\log \|b_I\|$}
        \put(52,35){\tiny $y=-1.x$}
        \put(26,16){\tiny $y=-1.x+1.5$}
      \end{overpic}
      \\
      (c)&(f)\\
    \end{tabular}
  \end{center}
  \caption{%
    \small
    Comparison of the numerical data between the $n=3$ case $F(\cE_T)$ for 
    $\cE_T=\{(1,0,0,1),(0,1,0,1),(0,0,1,1),(0,0,0,1)\}$ and the Sierpi\'nski tetrahedron $\cS$ about the asymptotic 
    behaviour of the bodies' volumes $V$ and surfaces $S$ sorted by their ``radii'' $\rho=V/S$ in descending order.
    No exact formulae are known for these plots so the lines shown above represent just an
    interpolation of the numerical data. 
    (a,c) Radii of the circles inscribed in the bodies vs k after sorting the radii in descending order
    for the $n=3$ case (left) and the Sierpi\'nski tetrahedron (right). (b,d) Areas of the bodies and (c,f) 
    their surface sorted according with their radii for the $n=3$ case (left) and the Sierpi\'nski tetrahedron (right).
    From these interpolation we obtain that for $F(\cE_T)$ we have $1.6\leq\dim_B F(\cE_T)\leq 2.9$ and for the 
    Sierpi\'nski tetrahedron we get the correct answer $\dim\cS=2$.
  }
\end{figure}

First we evaluate the Minkowsky dimension as the growth rate of the volume $V$ and surface $S$ 
of the bodies when sorted by the radius $\rho=V/S$. In this case we use the fact that~\cite{Gai06}
$$
\sum_{i=k_\epsilon+1}^\infty V_i\leq V_\epsilon \leq S\epsilon + H\epsilon^2+\epsilon\sum_{i=1}^{k_\epsilon} S_i + \sum_{i=k_\epsilon+1}^\infty V_i+\frac{4}{3}\epsilon^3
$$
where $V_\epsilon$ is the volume of the neighborhood of $F(\cE_T)$ of radius $\epsilon$, $k_\epsilon$
the integer such that $\rho_{k_\epsilon+1}\leq \epsilon \leq\rho_{k_\epsilon}$, $V_i$ and $S_i$ the volume
and surface of the body $Z_i$, $S$ and $H$ the surface and mean curvature of the starting tetrahedron.
From the numerical data (see fig.~\ref{fig:frdim-n3}(a-c)) we obtain that $\epsilon\asymp k^{-.57}$, 
$Ak^{-.77}\leq V_k \leq Bk^{-1.2}$ and $A'k^{-1.3}\leq S_k \leq B'k^{-.5}$, so that
$A''\epsilon^{1.34}\leq V_\epsilon \leq B''\epsilon^{.125}$ and therefore 
$1.66\leq \dim_MF\leq2.75$. Unfortunately, unlike in the $n=2$ case, we are not able to exclude
from this bounds that the fractal has integer dimension $\dim_M F=2$; note that this is exactly what
happens for the Tetrix, i.e. the three-dimensional analog of the Sierpinki triangle.

\begin{figure}
  \label{fig:labelsDim}
  \begin{center}
    \begin{tabular}{cc}
      \begin{overpic}[tics=5,width=5.3cm]{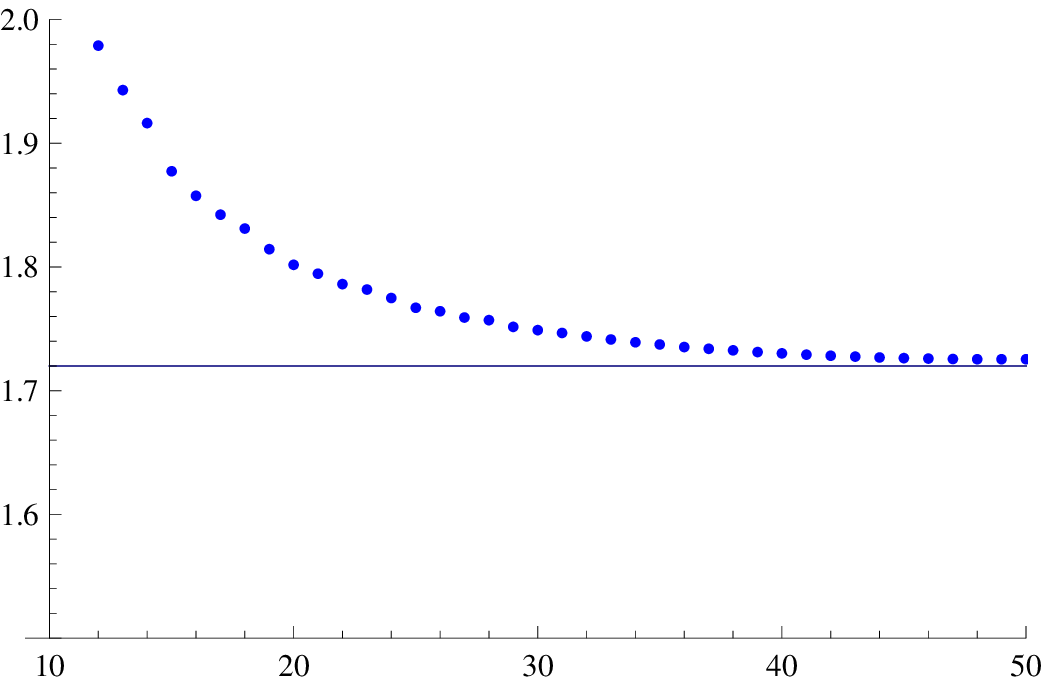}
        \put(19,55){\tiny $2-\frac{\log V_\epsilon}{\log\epsilon}$}
        \put(80,-3){\tiny $-\log\epsilon$}
      \end{overpic}
      &
      \begin{overpic}[tics=5,width=5.5cm]{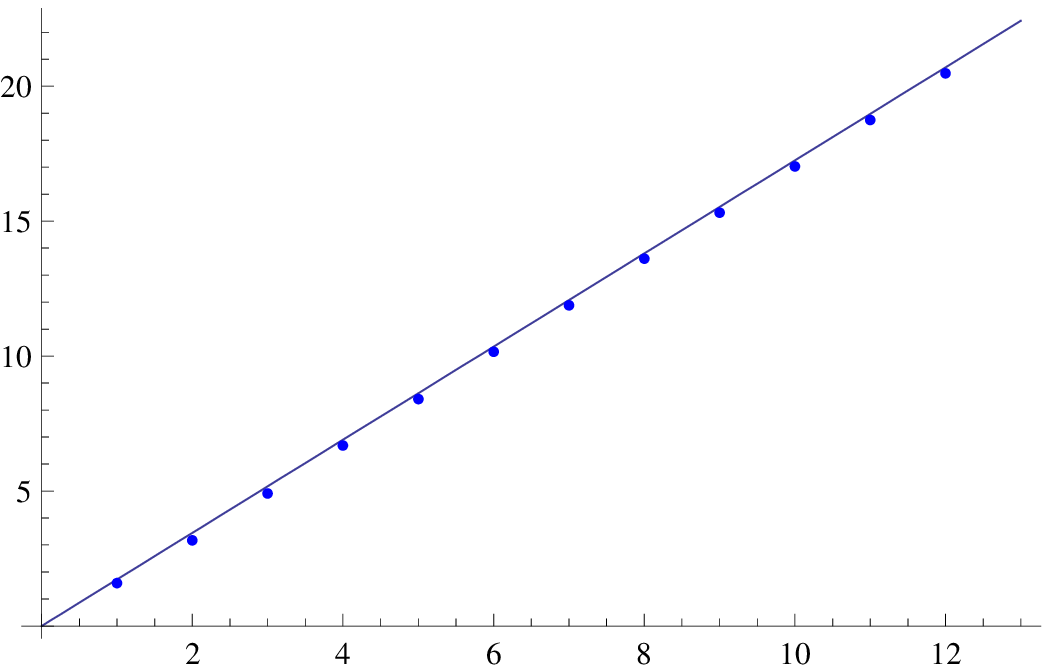}
        \put(9,55){\tiny $\log N_\epsilon$}
        \put(80,-3){\tiny $-\log\epsilon$}
        \put(40,40){\tiny $y=1.7x$}
      \end{overpic}
      \\
      (a)&(b)\\
      \begin{overpic}[tics=5,width=5.3cm]{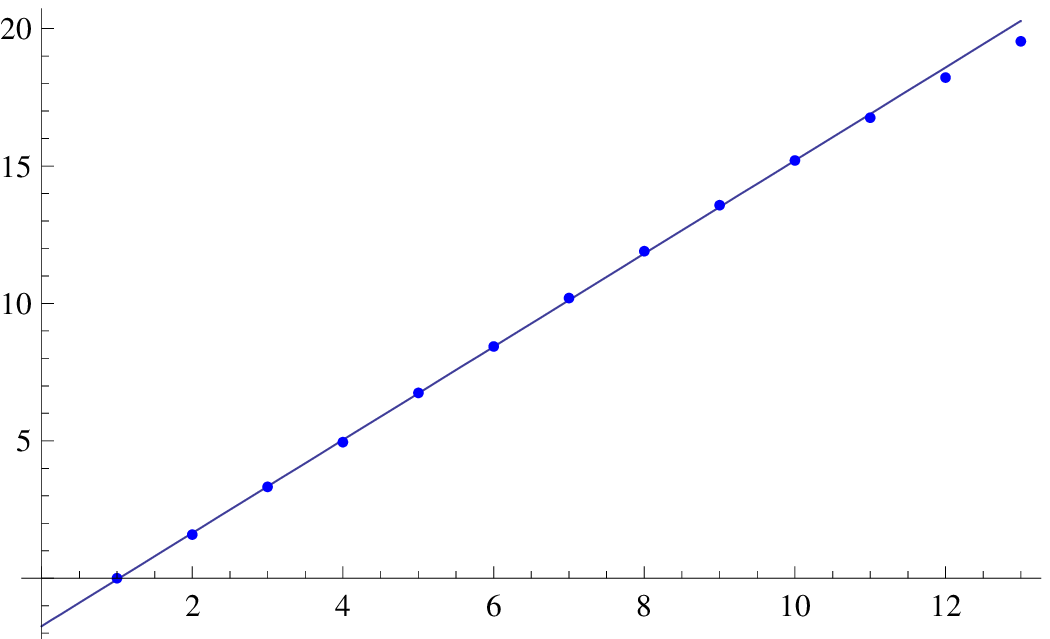}
        \put(9,55){\tiny $\log N_\epsilon$}
        \put(80,-3){\tiny $-\log\epsilon$}
        \put(40,40){\tiny $y=1.69x$}
      \end{overpic}
      &
      \begin{overpic}[tics=5,width=5.5cm]{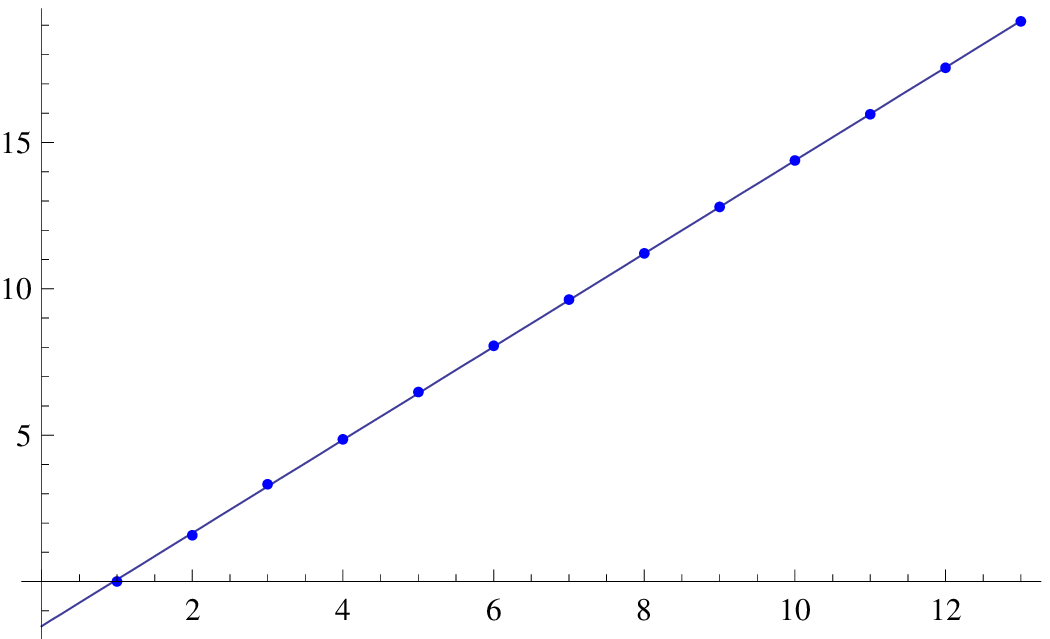}
        \put(9,55){\tiny $\log N_\epsilon$}
        \put(80,-3){\tiny $-\log\epsilon$}
        \put(40,40){\tiny $y=1.59x$}
      \end{overpic}
      \\
      (c)&(d)\\
      \begin{overpic}[tics=5,width=5.3cm]{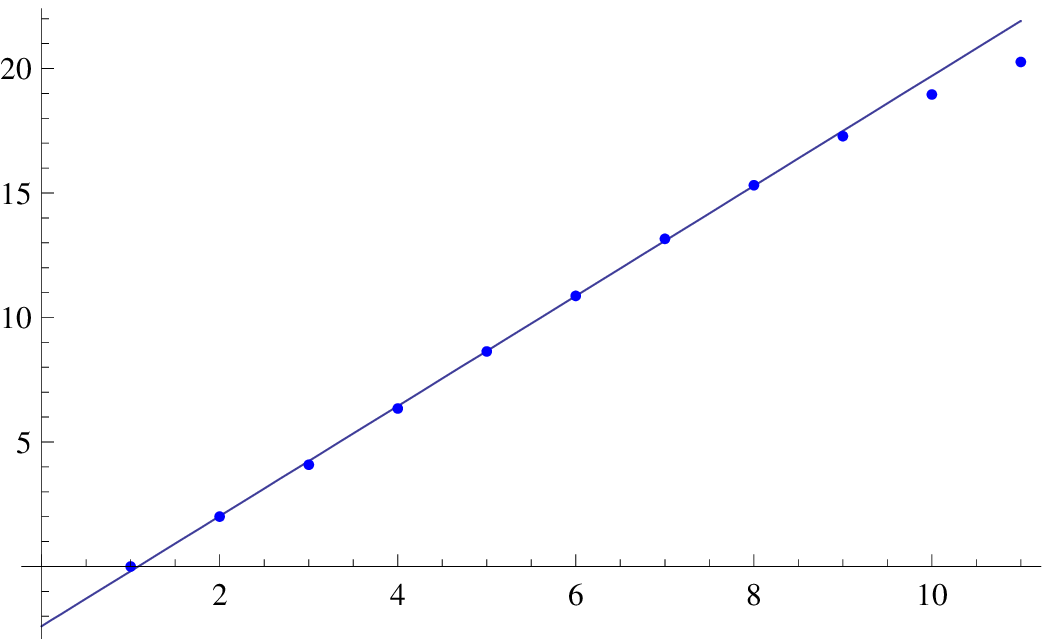}
        \put(9,55){\tiny $\log N_\epsilon$}
        \put(80,-3){\tiny $-\log\epsilon$}
        \put(40,40){\tiny $y=2.20x$}
      \end{overpic}
      &
      \begin{overpic}[tics=5,width=5.5cm]{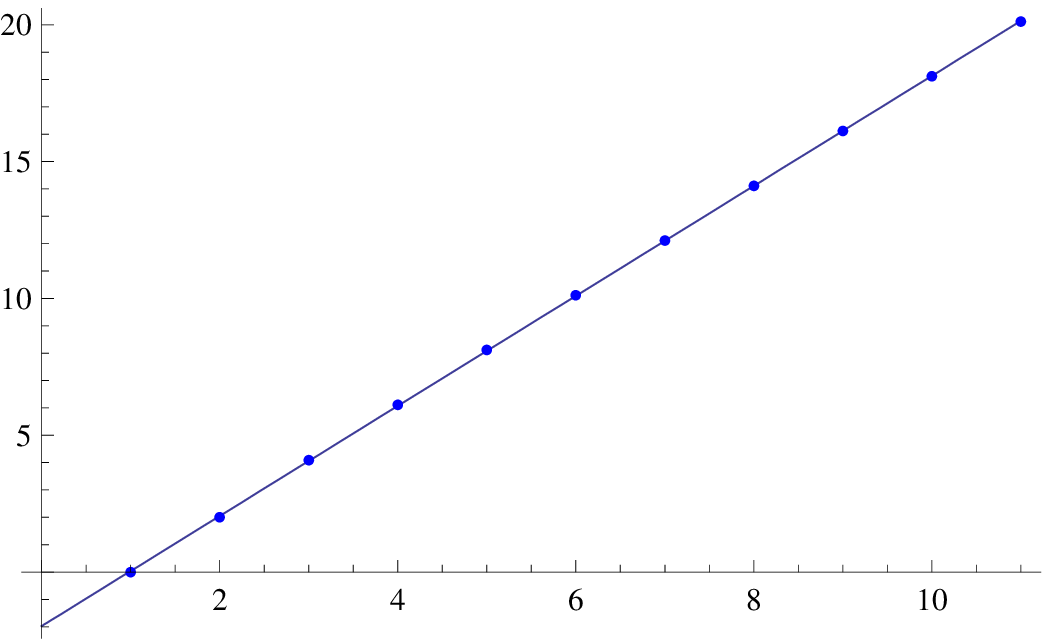}
        \put(9,55){\tiny $\log N_\epsilon$}
        \put(80,-3){\tiny $-\log\epsilon$}
        \put(40,40){\tiny $y=2.01x$}
      \end{overpic}
      \\
      (e)&(f)\\
    \end{tabular}
  \end{center}
  \caption{%
    \small
    Plots relative to the evaluation of fractal dimensions of $F_\cC=F(\cE_\cC)\subset\RPt$ 
    (see Section~\ref{sec:n2}) and $F_T=F(\cE_T)\subset\RPT$ (see Section~\ref{sec:n3}).
    (a) Evaluation of the Minkowsky dimension by direct numerical computation of the area $V_\epsilon$
    of the of the $\epsilon$ neighborhood of $F_\cC$, based on the fact that, for ``nice'' fractals $F$,
    $\dim F=\lim_{\epsilon\to0^+}\left[2-\frac{\log V_\epsilon}{\log\epsilon}\right]$.
    (b) Evaluation of the box-counting dimension $d_{bc}(F_\cC)$ by direct computation of the number 
    of squares $N_\epsilon$ needed to cover $F^{12}$, the 12-th order approximation of $F_\cC$, for 
    $\epsilon=2^{-k}$, $k=1,\cdots,12$. The data strongly suggest that $d_{bc}(F_\cC)\simeq1.7$.
    (c) Evaluation of the box-counting dimension $d'_{bc}(F_\cC)$ of the set of barycenters 
    (considered as points in $\RPt$) of the bodies $Z(\cE_I)$, $\cE_I\in T(\cE_\cC)$. 
    By Corollary~\ref{cor:bardensity}, the closure of this set is the union of $F(\cE_\cC)$ 
    with a 1-dim. set. We get, in excellent agreement with (b), $d'_{bc}(F_\cC)\simeq1.69$.
    (d) Same evaluation as in (c) in case of the Sierpinki triangle. The estimated dimension 
    $d'_{bc}\simeq1.59$ is in perfect agreement with the exact result $\log3/log2\simeq1.585$.
    (e) Evaluation of the box-counting dimension $d'_{bc}(F_T)$ of the set of barycenters 
    (considered as points in $\RPT$) of the bodies $Z(\cE_I)$, $\cE_I\in T(\cE_T)$. 
    By Corollary~\ref{cor:bardensity}, the closure of this set is the union of $F(\cE_T)$ 
    with a 2-dim. set. We get $d'_{bc}(F_T)\simeq2.20$.
    (f) Same evaluation as in (e) in case of the Tetrix (three-dim. analogue of the Sierpinki triangle). 
    The estimated dimension $d'_{bc}\simeq2.01$ is in excellent agreement with the exact result $2$.
  }
\end{figure}

Next we use Corollary~\ref{cor:bardensity} and evaluate numericlaly the box-counting dimension
$d_{bc}$ of the set $B=\{[b_I(\cE_T)]\}\subset\RPT$, namely the set of barycenters of the tree
$T(\cE_T)$ considered as points in the projective three-space. Analogously to the case $n=2$, 
the closure of $B$ is the union of $F(\cE_T)$ with the (two-dimensional) boundaries 
of the bodies, a dimension higher than two must be due to the points in $F(\cE_T)$.
We obtain $d_{bc}\simeq2.20$ (see Fig.~\ref{fig:labelsDim}). Notice that the very same
method, applied to the Tetrix, whose dimension is $d_T=2$, gives the quite close
result $d_T\simeq 2.01$.

We could not get useful information from the other two methods used in the $n=2$ case. 
In conclusion, the two numerical results we obtained are compatible with each other 
and the estimate of the box-counting dimension of barycenters is sufficiently far from
integer to make us think that, unlike the Tetrix, this fractal may have Hausdorff dimension 
higher than 2.

\begin{figure}
  \label{fig:asymp-n3}
  \begin{center}
    \begin{tabular}{cc}
      \begin{overpic}[tics=5,width=5.5cm]{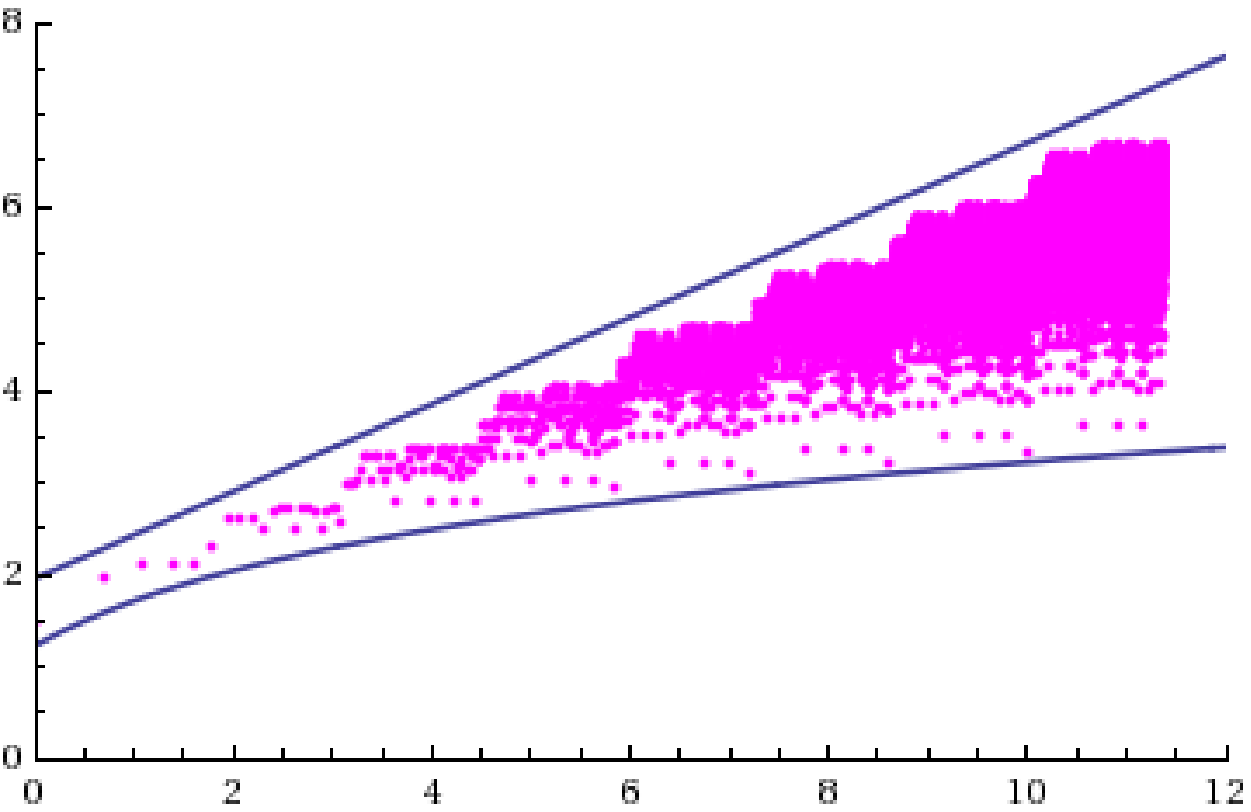}
        \put(6,55){\tiny $\log\|b_k\|$}
        \put(80,-3){\tiny $\log k$}
        \put(17,45){\tiny $y=\log_4\alpha_4 x+\log A$}
        \put(50,17){\tiny $y=\log(\frac{3x}{\log4}+\log_4108)$}
      \end{overpic}
      &
      \begin{overpic}[tics=5,width=5.5cm]{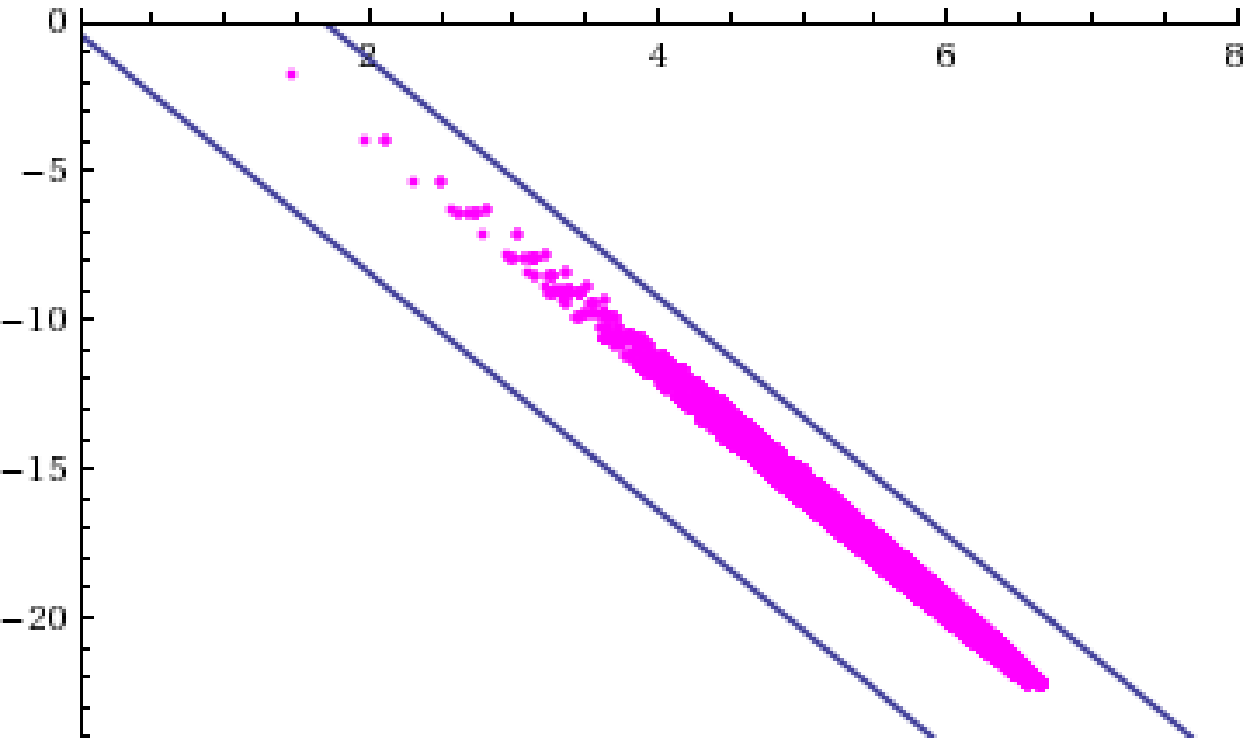}
        \put(9,37){\tiny $\log\mu(Z_I)$}
        \put(80,-3){\tiny $\log \|b_I\|$}
        \put(50,42){\tiny $y=-4x+\log(864)$}
        \put(17,16){\tiny $y=-4x-\log\frac{2}{3}$}
      \end{overpic}
      \\
      (a)&(b)\\
      \begin{overpic}[tics=5,width=5.5cm]{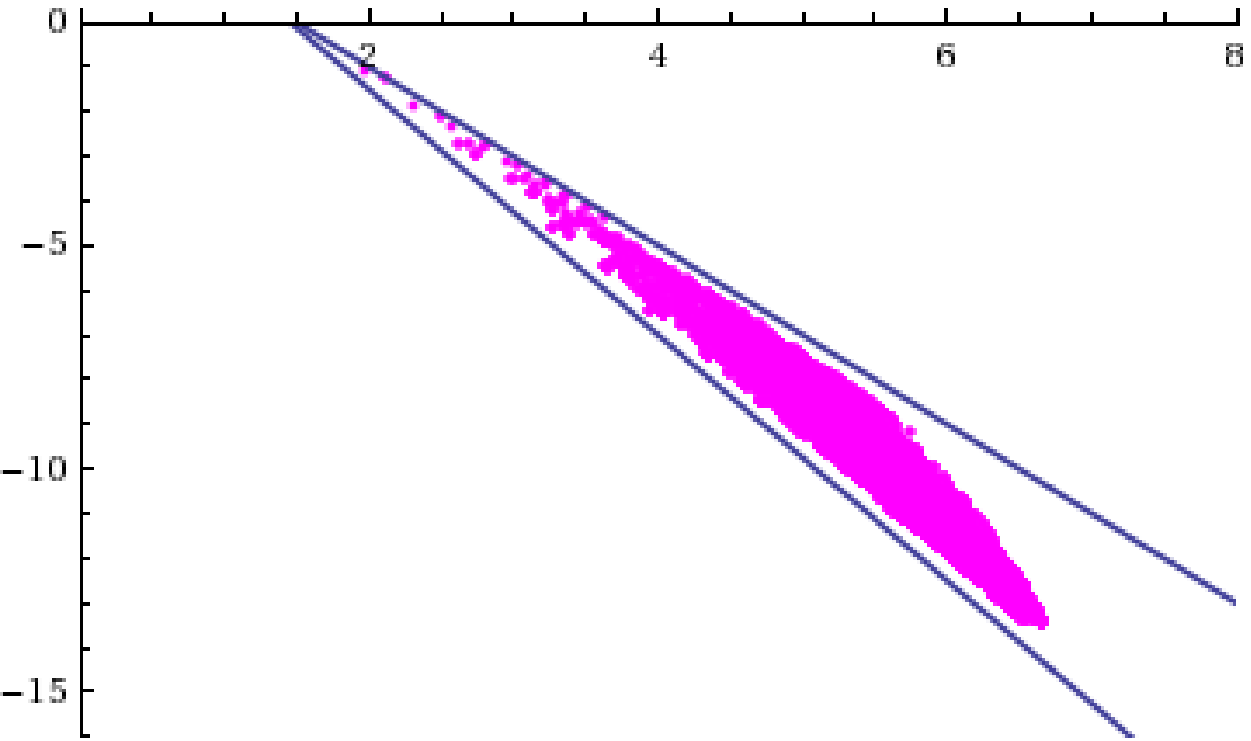}
        \put(9,37){\tiny $\log|Z_I|$}
        \put(80,-3){\tiny $\log \|b_I\|$}
        \put(55,42){\tiny $y=-x+2.2$}
        \put(25,18){\tiny $y=-1.33x+2.1$}
      \end{overpic}
      &
      \begin{overpic}[tics=5,width=5.5cm]{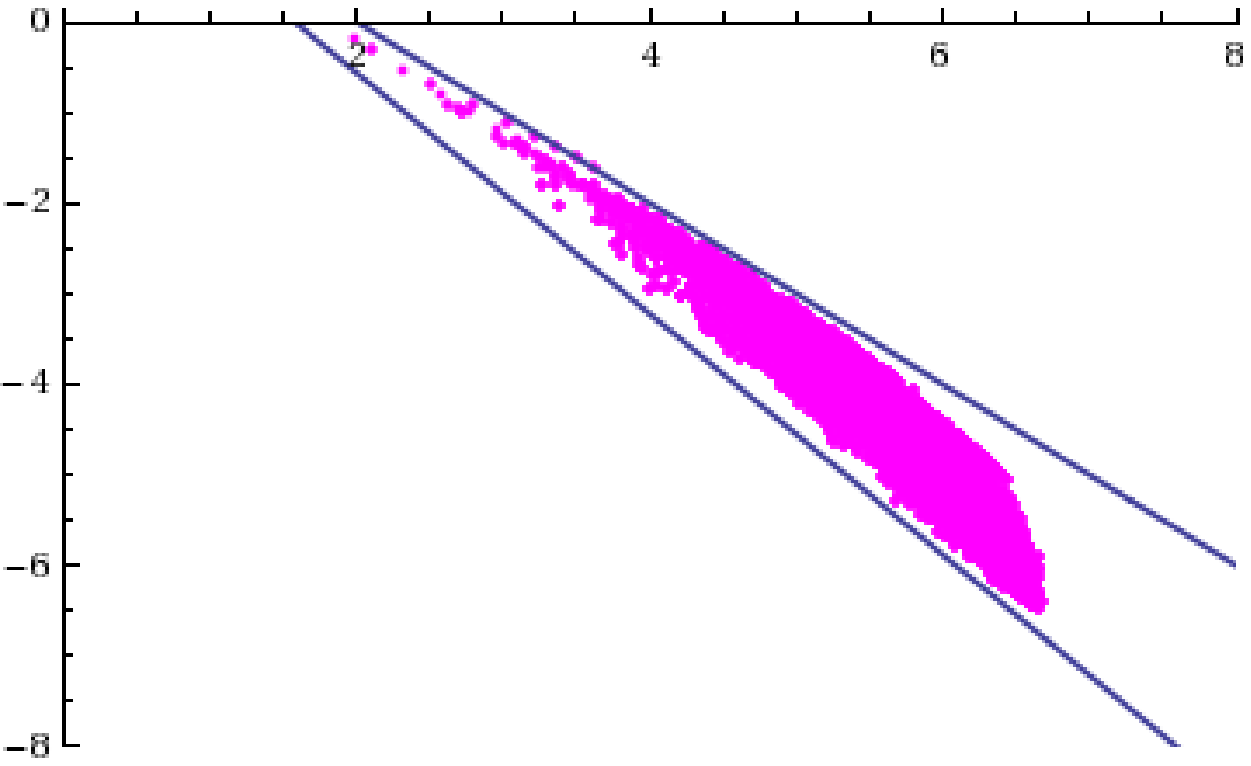}
        \put(9,37){\tiny $\log\mu_2(Z_I)$}
        \put(80,-3){\tiny $\log \|b_I\|$}
        \put(62,42){\tiny $y=-2.1x+3.5$}
        \put(30,23){\tiny $y=-2.8x+4$}
      \end{overpic}
      \\
      (c)&(d)\\
    \end{tabular}
  \end{center}
  \caption{%
    \small
    log-log plots for the main quantities in the $n=3$ case for 
    $\cE_T=\{(1,0,0,1),(0,1,0,1),(0,0,1,1),(0,0,0,1)\}$ 
    (a) Barycenters norms vs indices -- as explained in fig.~\ref{fig:plots-n1} the $b_k$ are ordered 
    according to the natural order induced by the tree, so that 
    $3k+4\leq\|b_\frac{4^k-1}{3}\|\leq2\alpha_3^4\frac{(1-\beta_4)(1-\bar\beta_4)(1-\gamma_4)}{(\alpha-\beta_4)(\alpha-\bar\beta_3)(\alpha-\gamma_4)}\alpha_4^k=A\alpha_4^k$
    and therefore
    $\frac{3}{\log4}\log k+\log_4108 \leq \|b_k\| \leq Ak^{\log_4\alpha_4}$.
    (b) Bodies' volumes vs barycenters norms -- the lines bounding the numerical data come immediately 
    from the inequalities in Theorem~\ref{thm:bodyVol-n2}.
    (c) Bodies' diameters vs barycenters norms and (d) bodies' surfaces areas; for these quantities
    no exact formulae are known so the lines shown represent just an interpolation of the numerical data. 
  }
\end{figure}
%
%
%
%
\section{Acknowledgments}
The author gladly thanks the IPST (www.ipst.umd.edu) and the Dept. of Mathematics of the UMD (USA) 
(www.math.umd.edu) for their hospitality in the Spring Semester 2007 and for financial support.
Numerical calculations were made on Linux PCs kindly provided by the UMD Mathematics Dept. 
and by the Cagliari section of INFN (www.ca.infn.it), that the author also thanks for financial support. 
The author finally warmly thanks S.P. Novikov and B. Hunt for several fruitful discussions during his 
stay at UMD and especially I.A. Dynnikov for many discussions on the subject, 
for helping clarifying the structure of the fractal when it was first discovered and for 
proofreading the manuscript.
\bibliography{fractal}
\end{document}